\input epsf

\documentclass[12pt]{article}
\usepackage{subfigure}
\usepackage[font=scriptsize]{caption}

\newcommand{\sign}{\operatorname{sign}}

\usepackage{latexsym,amsmath,amssymb,amsfonts,amscd,multirow}
\usepackage{epsfig,subfigure,verbatim,epstopdf,graphics}
\usepackage{color}

\newtheorem{thm}{Theorem}[section]

\newtheorem{exa}[thm]{Example}

\newcommand{\bx}{\mathbf{x}}

\newfont{\iams}{msbm9}

\newcommand{\commentbis}[1]{}
\newcommand{\be}{\begin{eqnarray}}
\newcommand{\ee}{\end{eqnarray}}
\newcommand{\beno}{\begin{eqnarray*}}
\newcommand{\eeno}{\end{eqnarray*}}
\newcommand{\barr}[1]{\begin{array}{#1}}
\newcommand{\earr}{\end{array}}

\newcommand{\beq}{\begin{equation}}
\newcommand{\eeq}{\end{equation}}
\newcommand{\beqa}{\begin{eqnarray}}
\newcommand{\eeqa}{\end{eqnarray}}

\newcommand{\bn}{{\bf n}}

\newcommand{\mT}{{\mathcal T}}

\newcommand{\jp}{{j+\frac{1}{2}}}
\newcommand{\jm}{{j-\frac{1}{2}}}
\newcommand{\ip}{{i+\frac{1}{2}}}
\newcommand{\im}{{i-\frac{1}{2}}}
\newcommand{\bt}{{\bf t }}

\newcounter{nsez}

\title
{  A New Discontinuous Galerkin Finite Element Method for
Directly Solving the 
Hamilton-Jacobi Equations\footnote{Research supported by  NSF grants DMS-1217563, DMS-1318186, AFOSR grant FA9550-12-1-0343 and the startup fund from Michigan State University.}}

\author{  Yingda Cheng
\thanks{Department of Mathematics, Michigan State University,
East Lansing, MI 48824 U.S.A.
 {\tt ycheng@math.msu.edu}}
  \and
Zixuan Wang
\thanks{Department of Mathematics, Michigan State University,
East Lansing, MI 48824 U.S.A.
 {\tt wangzix1@msu.edu}}
}

\date{\today}

\begin{document}

%======================================

\maketitle

\begin{abstract}
In this paper, we improve upon the discontinuous Galerkin (DG) method for Hamilton-Jacobi (HJ) equation with convex Hamiltonians in \cite{Cheng_07_JCP_HJ} and develop a new   DG method for directly solving the general HJ equations. 
The new method avoids the reconstruction of the solution across elements by utilizing the Roe speed at the cell interface.
Besides, we propose an entropy fix by adding penalty terms  proportional to the jump of the normal derivative of the numerical solution. The particular form of the entropy fix  was inspired by the  Harten and Hyman's entropy fix \cite{harten1983self} for Roe scheme for the conservation laws.    The resulting scheme is compact, simple to implement even on unstructured meshes, and is demonstrated to work for nonconvex Hamiltonians.
Benchmark numerical experiments in one dimension and two dimensions are provided to validate the performance of the method.

\end{abstract}


\section{Introduction}

In this paper, we consider the numerical solution of the time-dependent
Hamilton-Jacobi (HJ) equation
\begin{equation}
\label{hj}
\varphi_t+H(\nabla_\bx \varphi, \, \bx)=0,
 \quad \varphi(\bx, 0)=\varphi^0(\bx), \quad \bx \in \Omega \in \mathbb{R}^d
\end{equation}
with suitable boundary conditions on $\partial \Omega$. The HJ equation arises in many applications, e.g.,
optimal control, differential games, crystal growth, image
processing and calculus of variations. The solution of such equation
may develop discontinuous derivatives in finite time even when the
initial data is smooth. The viscosity solution \cite{Crandall_1983_TAMS_HamiJaco, crandall1984some} was introduced as the
unique physically relevant solution, and has been the focus of many
numerical methods. Starting from  \cite{Crandall_84_MC_HJ, souganidis1985approximation},    finite difference  methods such as
essentially non-oscillatory (ENO) \cite{osher1988fronts, Osher_1991_SIAM_NonOscill} or weighted ENO (WENO)  methods \cite{Jiang_1999_SIAM_WENO, Zhang_2003_SIAM_WENO}  have been developed to solve the HJ equation.  Those finite difference methods
work quite efficiently for Cartesian meshes, however they lose the advantage of simplicity on unstructured meshes \cite{abgrall1996numerical, Zhang_2003_SIAM_WENO}.

Alternatively, the Runge-Kutta discontinuous Galerkin (RKDG) method, originally devised
to solve the conservation laws \cite{Cockburn_2001_RK_DG}, is more flexible for
arbitrarily unstructured meshes. The first work of DG methods for HJ equations \cite{Hu_1999_SIAM_DG_FEM,
Li_2005_AML_DG_HamiJaco} relies on solving the conservation law system satisfied by the
derivatives of  the solution.  The methods  work well numerically even on unstructured mesh, with 
provable stability results for certain special cases, and were later generalized in eg. \cite{guo2011local, chen2007adaptive}.  Unfortunately, the procedure of recovering $\varphi$ from its derivatives has made the algorithm indirect and
complicated. In contrast, the design of  DG methods for directly solving the HJ equations is appealing but challenging, because the HJ equation is not written in the conservative form, for which the framework of  DG methods could easily apply. In \cite{Cheng_07_JCP_HJ},   a DG method for directly
solving the HJ equation   was developed.
This scheme has provable
stability and error estimates for linear equations and demonstrates
good convergence to the viscosity solutions for nonlinear equations.
However, in entropy violating cells, a correction based on the
schemes in \cite{Hu_1999_SIAM_DG_FEM, Li_2005_AML_DG_HamiJaco} is
necessary to guarantee stability of the method, and moreover, the method in \cite{Cheng_07_JCP_HJ} only works for equations with convex Hamiltonians. 
Later, this algorithm was applied to solve front propagation problems \cite{Cheng_09_front} and front propagation with obstacles \cite{Cheng_obs}, in which simplified implementations of the entropy fix procedure were proposed.
Meanwhile, central DG \cite{li2010central} and local DG \cite{yan2011local} methods were recently developed for the HJ equation. Numerical experiments demonstrate that both methods work for nonconvex Hamiltonians. In addition, the first order version of the local DG method \cite{yan2011local} reduces to the monotone schemes and thus has provable convergence properties. However, the central DG methods based on overlapping meshes are difficult to implement on unstructured meshes, and the local DG methods still need to resort to the information about the derivatives of $\varphi$, making the method less direct in computation. $L^2$ error estimates for smooth solutions of the DG method  \cite{Cheng_07_JCP_HJ} and local DG method  \cite{yan2011local} have been established in \cite{xiong2013priori}. For    recent developments of high order and DG methods for HJ equations, one can refer to the review papers \cite{shu2013survey,shu2007high}.

In this paper, we
improve upon the DG scheme in \cite{Cheng_07_JCP_HJ} and develop a new   DG method for directly solving the general HJ equation. Based on the observation that the method in \cite{Cheng_07_JCP_HJ} is closely related to Roe's linearization, we  use interfacial terms involving the Roe speed and develop a new entropy fix that was inspired by the  Harten and Hyman's entropy fix \cite{harten1983self} for Roe scheme for the conservation laws. 
The new method has the following advantages. Firstly, the scheme is shown to work on unstructured meshes even for nonconvex Hamiltonian. Secondly, the method is simple to implement. The cumbersome $L^2$ reconstruction of the solutions' derivative at the cell interface in \cite{Cheng_07_JCP_HJ} is avoided, and  the entropy fix is   automatically incorporated by the added jump terms in the derivatives of the numerical solution. Finally, the scheme is direct and compact, and the computation only needs the information about the current cell and its immediate neighbors.

The rest of the paper is  organized as follows: 
in Section \ref{sec:method}, we introduce the numerical schemes for one-dimensional and multi-dimensional HJ equations.   Section \ref{sec:numerical} is devoted to the discussion of the numerical results. We conclude and discuss about future work in Section \ref{sec:conclusion}.

%%%%%%%%%%NUMERICAL%%%%%%%%%%%%%
\section{Numerical Methods}
\label{sec:method}

In this section, we will describe the numerical methods. We follow the method of lines approach, and below we will only describe the semi-discrete DG schemes.   The resulting method of lines ODEs can be solved by the  total variation diminishing (TVD) 
 Runge-Kutta methods \cite{Shu_1988_JCP_NonOscill} or strong stability preserving Runge-Kutta methods \cite{Gottlieb_2001_SIAM_Stabi}.  
% We will start with the description for  one-dimensional equation first, and then introduce the methods in multi-dimensional case.
  
\subsection{Scheme in one dimension}
In this subsection, we will start by designing the scheme for one-dimensional HJ equation.
In this case, (\ref{hj}) becomes
\begin{equation}
\label{1dhj}
\varphi_t+H(\varphi_x,x)=0,\qquad
\varphi(x,0)=\varphi^0(x).
\end{equation}
Assume the computational domain is $[a,b]$,   we will
divide it into $N$ cells as follows
\begin{equation}
\label{cell1}
a=x_{\frac{1}{2}}<x_{\frac{3}{2}}< \ldots <x_{N+\frac{1}{2}}=b.
\end{equation}
Now  the cells and their centers are defined as
\begin{equation}
\label{cell2}
I_j=(x_{j-\frac{1}{2}},x_{j+\frac{1}{2}}), \qquad
x_j=\frac{1}{2} \left( x_{j-\frac{1}{2}}+x_{j+\frac{1}{2}} \right) , \quad j=1,\cdots, N
\end{equation}
and the mesh sizes are
\begin{equation}
\label{cell3}
\Delta x_j=x_{j+\frac{1}{2}}-x_{j-\frac{1}{2}}, \qquad
h =\max_j \Delta x_j.
\end{equation}
The  DG approximation space is
\begin{equation}
\label{space}
V_h^k=\{ \upsilon : \, \upsilon|_{I_j} \in P^k(I_j), \,\,
j=1,\cdots, N\}
\end{equation}
where $P^k(I_j)$ denotes all polynomials of degree at most $k$ on $I_j$, and we
let $H_1=\frac{\partial H}{\partial \varphi_x}$ be the partial derivative of the Hamiltonian with respect to $\varphi_x$. 

 To introduce the scheme, we need to define several quantities at the cell interface where the DG solution is discontinuous. If $x_*$ is a point located at the cell interface, then $\varphi_h \in V_h^k$ would be discontinuous at $x_*$. We can then define the Roe speed at $x_*$ to be
\[\tilde{H}_{\varphi_h} (x_*):=
 \left\{
  \begin{array}{ll}
\frac{H((\varphi_h)_x(x_*^+), x_*^+)-H((\varphi_h)_x(x_*^-), x_*^-)}{(\varphi_h)_x(x_*^+)- (\varphi_h)_x(x_*^-)},  & \textrm{if} \,(\varphi_h)_x(x_*^+) \neq  (\varphi_h)_x(x_*^-)\\
\frac{1}{2}(H_1((\varphi_h)_x(x_*^+), x_*^+)+H_1((\varphi_h)_x(x_*^-), x_*^-)), & \textrm{if} \,(\varphi_h)_x(x_*^+) =  (\varphi_h)_x(x_*^-). 
  \end{array} \right.\]
In the notations above, we use superscripts $+$, $-$ to denote the right, and left limits of a function.
Notice that in order for the above definition to make sense, we restrict our attention to $k \geq 1$ case, i.e. piecewise linear polynomials and above. 

Similar to Harten and Hyman's entropy fix \cite{harten1983self}, to detect entropy violating cell,  we introduce \beno
\delta_{\varphi_h}  (x_*) :=
\max \left (0,  \tilde{H}_{\varphi_h} (x_*)-H_1((\varphi_h)_x(x_*^-), x_*^-), H_1((\varphi_h)_x(x_*^+), x_*^+)- \tilde{H}_{\varphi_h} (x_*)\right)
\eeno
and
$$S_{\varphi_h} (x_*) :=\max \left (\delta_{\varphi_h}  (x_*), |\tilde{H}_{\varphi_h} (x_*)| \right).$$
We can see that $S_{\varphi_h} (x_*) \neq  |\tilde{H}_{\varphi_h} (x_*)|$ only if $ |\tilde{H}_{\varphi_h} (x_*)| < \delta_{\varphi_h}  (x_*)$.

Now we are ready to formulate our DG scheme for (\ref{1dhj}). We look for
$\varphi_h(x,t) \in V_h^k$, such that
\begin{align}
\label{1dscheme}
   & \int_{I_j} (\partial_t \varphi_h(x,t)+H(
        \partial_x \varphi_h(x,t),x))v_h(x) \, dx \nonumber \\
 &+   \min \left (\tilde{H}_{\varphi_h} (x_\jp), 0 \right)
    [\varphi_h]_\jp (v_h)_{j+\frac{1}{2}}^- \nonumber   + 
         \max \left (\tilde{H}_{\varphi_h} (x_\jm), 0 \right)
    [\varphi_h]_\jm(v_h)_{j-\frac{1}{2}}^+ \nonumber \\
    &- C \Delta x_j (S_{\varphi_h} (x_\jp)-|\tilde{H}_{\varphi_h} (x_\jp)| )[(\varphi_h)_x]_\jp (v_h)_{j+\frac{1}{2}}^-   \notag \\
        &-C \Delta x_j (S_{\varphi_h} (x_\jm)-|\tilde{H}_{\varphi_h} (x_\jm)| )[(\varphi_h)_x]_\jm (v_h)_{j-\frac{1}{2}}^+   \notag \\
 & = 0, \qquad \forall \,j=1,\ldots, N
\end{align}
holds for any $v_h \in V_h^k$ with $k \geq 1$. Here $[u]=u^+-u^-$ denotes the jump of a function at the cell interface, $\Delta x_j$ serves as the scaling factor.  $C$ is a positive  penalty parameter. The detailed discussion about the choice of $C$ is contained in Section  \ref{sec:numerical}. In particular, we find that $C=0.25$ works well in practice.

Next, we provide interpretation of the scheme \eqref{1dscheme} for the   linear HJ equation with variable coefficient
\begin{equation*}
\varphi_t+a(x)\varphi_{x}=0
\end{equation*}
to illustrate the main ideas. Firstly, if $a(x) \in C^1$, then 
$$\tilde{H}_{\varphi_h} (x_\jp)= a(x_\jp),\, \delta_{\varphi_h}  (x_\jp)=0, \, S_{\varphi_h} (x_\jp)=|a(x_\jp)| \quad \forall j=1,\ldots, N,$$ 
therefore scheme \eqref{1dscheme} reduces to

\begin{align}
\label{linearscheme}
  &\int_{I_j} ( \partial_t \varphi_h(x,t) + a(x)\partial_x \varphi_h(x,t) ) v_h(x) \,dx  \nonumber \\
  &+   \min \left (a(x_{j+\frac{1}{2}}), 0 \right)  [\varphi_h]_{j+\frac{1}{2}} (v_h)_{j+\frac{1}{2}}^- +\max \left (a(x_{j-\frac{1}{2}}), 0 \right)[\varphi_h]_{j-\frac{1}{2}} (v_h)_{j-\frac{1}{2}}^+  = 0, \quad j=1,\ldots, N.
\end{align}
This is the same DG method as in \cite{Cheng_07_JCP_HJ}, and it   is the standard DG schemes for the conservation laws with source term
$$\varphi_t+(a(x)\varphi)_x=a_x(x) \varphi$$
with Roe flux. Therefore, stability and error estimates could be established \cite{Cheng_07_JCP_HJ}.

On the other hand, when $a(x)$ is no longer smooth, especially when $a(x)$ contains discontinuity at cell interfaces. The scheme \eqref{linearscheme} will produce entropy violating shocks in the solutions' derivative \cite{Cheng_07_JCP_HJ}. In this case, the penalty terms in \eqref{1dscheme} come into play, and the added viscosity enables us to capture the viscosity solution as demonstrated in Section \ref{sec:numerical}. In particular, suppose $a(x)$ is discontinuous at $x_\jp$, then
$$\tilde{H}_{\varphi_h} (x_\jp)= \frac{[a(x)\varphi_h]_\jp}{[\varphi_h]_\jp}.$$
If the entropy condition is violated at $x_\jp$, i.e. $a(x_\jp^-)<0<a(x_\jp^+)$, then  
$$\delta_{\varphi_h}  (x_\jp)=\max \left (0,  \frac{[a(x)\varphi_h]_\jp}{[\varphi_h]_\jp}-a(x_\jp^-), a(x_\jp^+)- \frac{[a(x)\varphi_h]_\jp}{[\varphi_h]_\jp}\right)>0,$$
and $S_{\varphi_h}  (x_\jp)>0$. When $|\tilde{H}_{\varphi_h} (x_\jp)|< \delta_{\varphi_h}  (x_\jp)$, the penalty term $$-C \Delta x_j (S_{\varphi_h} (x_\jp)-|\tilde{H}_{\varphi_h} (x_\jp)| )[(\varphi_h)_x]_\jp (v_h)_{j+\frac{1}{2}}^-  $$ will be nonzero.  On the other hand, if $a(x_\jp^-)>0>a(x_\jp^+)$, corresponding to a shock in $\varphi_x$, we know the Roe scheme \eqref{linearscheme} could correctly capture this behavior. In fact, in this case
$$\delta_{\varphi_h}  (x_\jp)=\max \left (0,  \frac{[a(x)\varphi_h]_\jp}{[\varphi_h]_\jp}-a(x_\jp^-), a(x_\jp^+)- \frac{[a(x)\varphi_h]_\jp}{[\varphi_h]_\jp}\right)=0,$$
and $S_{\varphi_h} (x_\jp)-|\tilde{H}_{\varphi_h} (x_\jp)|=0,$ the method \eqref{1dscheme} will reduce to \eqref{linearscheme}. 
Similar arguments extend to the nonlinear case for sonic expanding cells for convex Hamiltonians,
\begin{eqnarray*}
H_1(\varphi_x^-(x_{j+\frac{1}{2}})) < 0 < H_1(\varphi_x^+(x_{j+\frac{1}{2}})).
\end{eqnarray*}
The penalty term in \eqref{1dscheme} would be turned on automatically. 

Finally we remark that the key differences of the scheme \eqref{1dscheme} compared to the method in \cite{Cheng_07_JCP_HJ} are: (1) the $L^2$ reconstruction of $\varphi_h$ across two elements is avoided, and we use the Roe speed which could be easily computed. This is advantageous especially for multi-dimensional problems on unstructured meshes, as illustrated in Sections \ref{sec:rectangle} and \ref{sec:unstructured}. (2) The added penalty terms automatically treat the sonic points, and the key idea is to add the viscosity based on the jump in $(\varphi_h)_x$, but not $\varphi_h$ itself. This is natural considering the formation of  monotone schemes such as the Lax-Friedrichs scheme for HJ equation. We will verify in Section \ref{sec:numerical} that the penalty terms do not degrade the order of the accuracy of the numerical scheme.

%Finally, as noted in \cite{Cheng_07_JCP_HJ}, for nonlinear equations, and $k=0$, the method is not consistent. Also we can see that in the $P^0$ case, $(\varphi_h)_x$ is always zero, and it is not consistent for the nonlinear equation. Therefore, in this paper, we will only consider the case of $k \ge 1$.

%Motivated by Roe fix of Harten and Hyman for the conservation law, we introduced the penalty term to add in more viscosity to fix the entropy violation
%$$ C (S ((\varphi_h)_x^-, (\varphi_h)_x^+, x)-|\tilde{H} ((\varphi_h)_x^-, (\varphi_h)_x^+, x)| )[(\varphi_h)_x] (v_h)_{j+\frac{1}{2}}^- \Delta x_j  $$ and
% $$C (S ((\varphi_h)_x^-, (\varphi_h)_x^+, x)-|\tilde{H} ((\varphi_h)_x^-, (\varphi_h)_x^+, x)| )[(\varphi_h)_x] (v_h)_{j-\frac{1}{2}}^+ \Delta x_j  $$\\
%Remark:\\
%$S ((\varphi_h)_x^-, (\varphi_h)_x^+, x) =\max \left (\delta, |\tilde{H} ((\varphi_h)_x^-, (\varphi_h)_x^+, x)| \right) $ is not equal to $ |\tilde{H} ((\varphi_h)_x^-, (\varphi_h)_x^+, x)|$ only in the sonic expanding cells. So the penalty term is automatically turned on when the entropy fix is needed, otherwise it is silent.

%
%\textcolor{red}{TODO: how to choose $C$ optimally? When $C=0$, there is no entropy fix, and will produce entropy violating solution. If $C<0$, the scheme should blow up. Need to test the value of $C$ numerically, e.g. $C=0.25$ seems to be the old entropy fix, could try $C=0.1$ etc.}
%
%
%\textcolor{red}{TODO: need to fix the phase point for non convex Hamiltonians, like those in \cite{serna2009characteristic}}

\subsection{Scheme on two-dimensional Cartesian meshes}
\label{sec:rectangle}
In this subsection, we generalize the scheme to compute on two-dimensional Cartesian meshes.
Now   equation (\ref{hj}) is written as
\begin{equation}
\label{2dhj}
\varphi_t+H(\varphi_x,\varphi_y,x,y)=0, \quad \varphi(x,y,0)=\varphi^0(x,y), \quad (x, y) \in [a, b] \times [c, d].
\end{equation}

The rectangular mesh is defined by
\begin{equation}
\label{2dcell1}
a=x_{\frac{1}{2}}<x_{\frac{3}{2}}< \ldots <x_{N_x+\frac{1}{2}}=b ,
\qquad
c=y_{\frac{1}{2}}<y_{\frac{3}{2}}< \ldots <y_{N_y+\frac{1}{2}}=d
\end{equation}
and
\begin{eqnarray}
\label{2dcell2}
 I_{i,j}=[x_{i-\frac{1}{2}},x_{i+\frac{1}{2}}] \times 
[y_{j-\frac{1}{2}},y_{j+\frac{1}{2}}] ,\quad J_i=[x_{i-1/2},x_{i+1/2}], 
\quad K_j=[y_{j-1/2},y_{j+1/2}] &&\nonumber \\
 i=1,\ldots N_x, \quad j=1,\ldots N_y.&&
\end{eqnarray}
Let $$\Delta x_i=x_{i+1/2}-x_{i-1/2}, \quad  \Delta y_j=y_{j+1/2}-y_{j-1/2}, \quad h=\max( \max_i \Delta x_i, \max_j \Delta y_j).$$
We define the approximation space as
\begin{equation}
\label{space2}
V_h^k=\{ \upsilon : \, \upsilon|_{I_{i,j}} \in P^k(I_{i,j}), \,\,
i=1,\ldots N_x, \quad j=1,\ldots N_y\}
\end{equation}
where $P^k(I_{i,j})$ denotes all polynomials of degree at most $k$ on 
$I_{i,j}$ with $k \geq 1$.

Let us denote $H_1=\frac{\partial H}{\partial \varphi_x}$
and $H_2=\frac{\partial H}{\partial \varphi_y}$. Similar to the one-dimensional case, we need to introduce several numerical quantities at the cell interface.

If $x_*$ is   located at the cell interface in the $x$ direction, then $\varphi_h \in V_h^k$ is discontinuous at $(x_*, y)$ for any $y$, and we define the Roe speed in the $x$ direction at $(x_*, y)$ to be
\beno
&&\tilde{H}_{1, \varphi_h} (x_*, y):=\\
 &&\hspace{-0.3in}\left\{
  \begin{array}{ll}
\frac{H((\varphi_h)_x(x_*^+,y), \overline{(\varphi_h)_y}, x_*^+,y)-H((\varphi_h)_x(x_*^-,y), \overline{(\varphi_h)_y}, x_*^-, y)}{(\varphi_h)_x(x_*^+,y)- (\varphi_h)_x(x_*^-,y)},  & \textrm{if} \,(\varphi_h)_x(x_*^+,y) \neq  (\varphi_h)_x(x_*^-,y)\\
\frac{1}{2}(H_1((\varphi_h)_x(x_*^+,y), \overline{(\varphi_h)_y}, x_*^+, y)+H_1((\varphi_h)_x(x_*^-,y), \overline{(\varphi_h)_y}, x_*^-, y)), & \textrm{if} \,(\varphi_h)_x(x_*^+,y) =  (\varphi_h)_x(x_*^-,y),
  \end{array} \right. \notag
  \eeno
  where $$
\overline{(\varphi_h)_y}=\frac{1}{2} \left(
( \varphi_h)_y(x_*^+,y) + (\varphi_h)_y(x_*^-,y) \right) $$
is the average of the tangential derivative. Again,  we define 
\beno
&&\delta_{1, \varphi_h}  (x_*, y) :=\\
&&\hspace{-0.3in}\max \left (0,  \tilde{H}_{1,\varphi_h} (x_*,y)-H_1((\varphi_h)_x(x_*^-, y), \overline{(\varphi_h)_y}, x_*^-, y), H_1((\varphi_h)_x(x_*^+, y), \overline{(\varphi_h)_y}, x_*^+, y)- \tilde{H}_{1, \varphi_h} (x_*, y)\right)
\eeno
and
$$S_{1, \varphi_h} (x_*, y) :=\max \left (\delta_{1, \varphi_h}  (x_*, y), |\tilde{H}_{1, \varphi_h} (x_*,y)| \right).$$
Similarly, for $y_*$   located at the cell interface in the $y$ direction,   $\varphi_h \in V_h^k$ is discontinuous at $(x, y_*)$ for any $x$, and we define the Roe speed in the $y$ direction at $(x, y_*)$ to be
\beno
&&\tilde{H}_{2, \varphi_h} (x, y_*):=\\
 &&\hspace{-0.3in}\left\{
  \begin{array}{ll}
\frac{H(\overline{(\varphi_h)_x},( \varphi_h)_y(x,y_*^+) ,  x, y_*^+)-H(\overline{(\varphi_h)_x}, (\varphi_h)_y(x, y_*^-),   x, y_*^-)}{(\varphi_h)_y(x,y_*^+)- (\varphi_h)_y(x,y_*^-)},  & \textrm{if} \,(\varphi_h)_y(x, y_*^+) \neq  (\varphi_h)_y(x, y_*^-)\\
\frac{1}{2}(H_1(\overline{(\varphi_h)_x}, (\varphi_h)_y(x,y_*^+), x, y_*^+)+H_1(\overline{(\varphi_h)_x}, (\varphi_h)_y(x,y_*^-), x, y_*^-)), & \textrm{if} \,(\varphi_h)_y(x, y_*^+) =  (\varphi_h)_y(x, y_*^-), 
  \end{array} \right. \notag
  \eeno
  where $$
\overline{(\varphi_h)_x}=\frac{1}{2} \left(
( \varphi_h)_x(x,y_*^+) + (\varphi_h)_x(x, y_*^-) \right) $$
is the average of the tangential derivative. Again,  we define 
\beno
&&\delta_{2, \varphi_h}  (x, y_*) :=\\
&&\hspace{-0.3in}\max \left (0,  \tilde{H}_{2,\varphi_h} (x,y_*)-H_2(\overline{(\varphi_h)_x},(\varphi_h)_y(x, y_*^-),  x, y_*^-), H_2(\overline{(\varphi_h)_x},(\varphi_h)_y(x, y_*^+),   x, y_*^+)- \tilde{H}_{2, \varphi_h} (x, y_*)\right)
\eeno
and
$$S_{2, \varphi_h} (x, y_*) :=\max \left (\delta_{2, \varphi_h}  (x, y_*), |\tilde{H}_{2, \varphi_h} (x,y_*)| \right).$$

We now formulate our scheme as: find $\varphi_h(x,t) \in V_h^k$, such that
\begin{eqnarray}
\label{2dscheme}
    &&\int_{I_{i,j}} (\partial_t \varphi_h(x,y,t)+H(
        \partial_x \varphi_h(x,y,t),\partial_y 
\varphi_h(x,y,t),x,y))v_h(x,y)dxdy \nonumber \nonumber \\
&&+ 
       \int_{K_j} \min \left (\tilde{H}_{1, \varphi_h}(x_\ip, y), 0 \right) [\varphi_h](x_{i+\frac{1}{2}},y) v_h(x^-_{i+\frac{1}{2}},y) 
dy\nonumber \\
&&+ 
       \int_{K_j} \max \left (\tilde{H}_{1, \varphi_h}(x_\im, y), 0 \right) [\varphi_h](x_{i-\frac{1}{2}},y) v_h(x^+_{i-\frac{1}{2}},y) 
dy\nonumber \\
&&+  
        \int_{J_i} \min \left (\tilde{H}_{2, \varphi_h}(x, y_\jp), 0 \right) [\varphi_h](x, y_\jp) v_h(x,y_\jp^-) 
dx\nonumber \\
 && + 
        \int_{J_i} \max \left (\tilde{H}_{2, \varphi_h}(x, y_\jm), 0 \right) [\varphi_h](x, y_\jm) v_h(x,y_\jm^+) 
dx  \\
&&- C \Delta x_i
       \int_{K_j}\left (S_{1, \varphi_h}(x_\ip, y)-|\tilde{H}_{1, \varphi_h} (x_\ip, y)| \right) [(\varphi_h)_x](x_{i+\frac{1}{2}},y) v_h(x^-_{i+\frac{1}{2}},y) 
dy \nonumber \\
&&- C \Delta x_i
       \int_{K_j}\left (S_{1, \varphi_h}(x_\im, y)-|\tilde{H}_{1, \varphi_h} (x_\im, y)| \right) [(\varphi_h)_x](x_{i-\frac{1}{2}},y) v_h(x^+_{i-\frac{1}{2}},y) 
dy\nonumber \\
&&- C \Delta y_j
        \int_{J_i} \left (S_{2, \varphi_h}(x, y_\jp)-|\tilde{H}_{2, \varphi_h} (x, y_\jp)| \right) [(\varphi_h)_y](x, y_\jp) v_h(x,y_\jp^-) 
dx\nonumber \\
 && - C \Delta y_j
        \int_{J_i} \left (S_{2, \varphi_h}(x, y_\jm)-|\tilde{H}_{2, \varphi_h} (x, y_\jm)| \right) [(\varphi_h)_y](x, y_\jm) v_h(x,y_\jm^+) 
dx
=0 \notag
\end{eqnarray}
holds for any $v_h \in V_h^k$ with $k \geq 1$.  In practice, the volume and line integrals in \eqref{2dscheme} can be evaluated by Gauss quadrature formulas.
The main idea in \eqref{2dscheme} is that in the normal direction of the interface, we apply the ideas from the  one-dimensional case, but  tangential to the 
interface, we   use the average of the tangential derivatives from the two 
neighboring cells. 

\subsection{Scheme on general unstructured meshes}
\label{sec:unstructured}
In this subsection, we describe the scheme on general unstructured meshes for \eqref{hj}. Let $\mathcal{T}_h=\{K\}$ be a partition  of $\Omega$,  with $K$ being simplices.
We define the piecewise polynomial space
$$
V_h^{k}=\left\{v\in L^2(\Omega): v|_K\in P^k(K), \,\forall K \in\mT_h\right\},
$$
where $P^k(K)$ denotes the set of polynomials of  total degree at most $k$ on $K$ with $k \geq 1$. For any element $K$, and  edge in $\partial K$, we define $\bn_K$ to be the outward unit normal to the boundary of $K$, and $\bt_K$ is the unit tangential vector such that $\bn_K \cdot \bt_K=0$. In higher dimensions, i.e. $d>2$, $d-1$ tangential vectors need to  be defined.  
In addition, for any
 function $u \in V_h^k$, and $\bx \in \partial K$, we define the traces of $u_h$ from  outside and inside of the element $K$ to be
 $$u_h^{\pm}(\bx)=\lim_{\epsilon \downarrow 0} u_h(\bx\pm\epsilon \bn_k),$$
 and $[u_h](\bx)=u_h^+(\bx)-u_h^-(\bx)$, $\overline{u_h}(\bx)=\frac{1}{2}(u_h^+(\bx)+u_h^-(\bx))$. We also let $H_{\bn_K}=\nabla_{\nabla \varphi} H \cdot \bn_K$.

Now following the Cartesian case, we define, for any $\bx \in \partial K$,
\beno
&&\tilde{H}_{\bn_K, \varphi_h} (\bx):=\\
 &&\hspace{-0.5in}\left\{
  \begin{array}{ll}
\frac{H ((\nabla_\bx \varphi_h \cdot \bn_K)^+, \overline{\nabla_\bx \varphi_h \cdot \bt_K}, \bx^+)-H ((\nabla_\bx \varphi_h \cdot \bn_K)^-, \overline{\nabla_\bx \varphi_h \cdot \bt_K}, \bx^-)}{[\nabla_\bx \varphi_h \cdot \bn_K] (\bx) },  & \textrm{if} \,[\nabla_\bx \varphi_h \cdot \bn_K] (\bx) \neq  0\\
\frac{1}{2}(H_{\bn_K}((\nabla_\bx \varphi_h \cdot \bn_K)^+, \overline{\nabla_\bx \varphi_h \cdot \bt_K}, \bx^+)+H_{\bn_K}((\nabla_\bx \varphi_h \cdot \bn_K)^-, \overline{\nabla_\bx \varphi_h \cdot \bt_K}, \bx^-)), & \textrm{otherwise}, 
  \end{array} \right. \notag
  \eeno
 and
\beno
&&\delta_{\bn_K, \varphi_h}  (\bx) :=\\
&&\hspace{-0.8in}\max \left (0, \tilde{H}_{\bn_K, \varphi_h} (\bx)-H_{\bn_K}((\nabla_\bx \varphi_h \cdot \bn_K)^-, \overline{\nabla_\bx \varphi_h \cdot \bt_K}, \bx^-), H_{\bn_K}((\nabla_\bx \varphi_h \cdot \bn_K)^+, \overline{\nabla_\bx \varphi_h \cdot \bt_K}, \bx^+)- \tilde{H}_{\bn_K, \varphi_h} (\bx)\right),
\eeno

$$S_{\bn_K, \varphi_h} (\bx) :=\max \left (\delta_{\bn_K, \varphi_h}  (\bx), |\tilde{H}_{\bn_K, \varphi_h} (\bx)| \right).$$

Then we look for $\varphi_h \in V_h^k$, such that for each $K$,
\begin{eqnarray*}
\label{2dgeneral.scheme}
    \int_{K} ((\varphi_h)_t+H(\nabla_\bx \varphi_h, \, \bx))\,v_h \,d\bx 
 +    \int_{\partial K} \min \left ( \tilde{H}_{\bn_K, \varphi_h} (\bx) , 0 \right)[\varphi_h](\bx) v_h^-(\bx) ds\\
  -C \frac{\Delta K}{\Delta S_K}
       \int_{\partial K} ( S_{\bn_K, \varphi_h} (\bx)  -|\tilde{H}_{\bn_K, \varphi_h} (\bx)|)[\nabla_\bx \varphi_h \cdot \bn_K](\bx) v_h^-(\bx) ds   =0  
\end{eqnarray*}
for any test function $v_h \in V_h^k$ with $k \geq 1$, where $\Delta K$, $\Delta S_K$ are size of the element $K$ and edge $S_K$ respectively. 
In practice, the volume and edge integrals need to be evaluated by quadrature rules with enough accuracy. For example, we use quadrature rules that are exact for $(2k)$-th order polynomial for the volume integral, and quadrature rules that are exact for $(2k+1)$-th order polynomial for the edge integrals.

\section{Numerical Results}
\label{sec:numerical}
In this section, we provide numerical   results to demonstrate the high order accuracy and reliability of our schemes. In all numerical experiments, we use the third order TVD-RK method as the temporal discretization \cite{Shu_1988_JCP_NonOscill}.

\subsection{One-dimensional results}

In this subsection, we provide computational results for one-dimensional HJ equations.
\begin{exa}\rm We solve the following linear problem with smooth variable coefficient
\label{linsmth}
\begin{flalign}
\ \ \ \ \ &\begin{cases}
    \varphi_t+\sin(x)\varphi_{x}=0 \\
    \varphi(x,0)=\sin(x)\\
    \varphi(0,t)=\varphi(2\pi,t)
  \end{cases}&
\end{flalign}
Since $a(x)$ is smooth in this example,  the penalty term automatically vanishes and the choice of $C$ does not have an effect on the solution. We provide the numerical results for $P^1$,$P^2$ and $P^3$ polynomials  in Table \ref{table:linsmth}. The CFL numbers used are also listed in this table. For $P^3$ polynomials, we set $\Delta t=O(\Delta x^{\frac{4}{3}})$ in order to get comparable numerical errors in time as in space. From the results, we could clearly observe the optimal $(k+1)$-th order accuracy for $P^k$ polynomials.

% Table 3.1
\begin{table}[ht]
\caption{Errors and numerical orders of accuracy for Example \ref{linsmth} when using $P^k$, $k=1,2,3$, polynomials and third order Runge-Kutta 
 time discretization on a uniform mesh of $N$ cells. 
Final time $t=1$.}
\vspace{2 mm}
\centering
\begin{tabular}{c c c c c c c}
\hline
 &      & $ P^1$   &    $CFL=0.3$  & &  &   \\
\hline
N   & $L^1$ error & order & $L^2$ error & order & $L^\infty$ error & order  \\
\hline
40  & 1.20E-03 &      & 2.55E-03 &       & 1.52E-02 &   \\
80  & 3.07E-04 & 1.96 & 6.83E-04 &  1.90 & 4.32E-03 &  1.81  \\
160 & 7.84E-05 & 1.97 & 1.78E-04 &  1.94 & 1.14E-03 &  1.92 \\
320 & 1.99E-05 & 1.98 & 4.56E-05 &  1.97 & 2.94E-04 &  1.96  \\
640 & 5.03E-06 & 1.99 & 1.15E-05 &  1.98 & 7.43E-05 &  1.98  \\
\hline
 &     & $P^2$  & $CFL=0.1$  &  &  &   \\
 \hline
 40  & 4.76E-05 &      & 9.97E-05 &       & 5.23E-04 &        \\
80  & 5.97E-06 & 2.99 & 1.36E-05 &  2.88 & 8.77E-05 &  2.58  \\
160 & 7.48E-07 & 3.00 & 1.82E-06 &  2.90 & 1.35E-05 &  2.70  \\
320 & 9.38E-08 & 2.99 & 2.38E-07 &  2.93 & 1.96E-06 &  2.78  \\
640 & 1.18E-08 & 2.99 & 3.08E-08 &  2.95 & 2.72E-07 &  2.85  \\
\hline
 &  &    $P^3$ &    $CFL=0.05$ &  &  &   \\
 \hline
40  &  2.12E-06 &       & 5.13E-06 &      & 2.89E-05 &   \\
80  &  1.36E-07 &   3.97 & 3.49E-07 & 3.89 & 2.16E-06 &   3.75  \\
160 &  8.71E-09 &  3.97 & 2.30E-08 & 3.93 & 1.57E-07 &  3.79 \\
320 &  5.14E-10  & 4.09 & 1.35E-09 & 4.10 &  9.47E-09 &4.06\\
640 &  4.83E-12 &  6.75 & 9.06E-12 & 7.24 & 4.52E-11 &  7.73\\
1280 &  2.03E-13  &  4.58  & 2.96E-13 & 4.94&  1.42E-12 &5.00\\
\hline
\end{tabular}
\label{table:linsmth}
\end{table}

\end{exa}

\begin{exa}\rm \label{linnonsmth} We solve the following linear problem with nonsmooth variable coefficient \cite{Cheng_07_JCP_HJ}
\begin{flalign}
\ \ \ \ \ & \begin{cases}
    \varphi_t+\text{sign}(\cos(x))\varphi_{x}=0 \\
    \varphi(x,0)=\sin(x)\\
    \varphi(0,t)=\varphi(2\pi,t)
  \end{cases}&
\end{flalign}

\rm The viscosity solution of this example has a shock forming in $\varphi_x$ at $x=\pi / 2$, and a rarefaction wave at $x=3\pi / 2$. 
\rm We use this example  to demonstrate the effect of the choice of $C$ on the numerical solution. The solutions obtained with different values of $C$ are provided in Figure \ref{fig linnonsmth}. If we take the penalty constant $C=0$, that is, without entropy correction, the entropy condition is violated at the two cells neighboring $x=3\pi / 2$, and the numerical solution is not convergent towards the exact solution. As we slowly increasing the value of $C$, we could observe better and better convergence property. In particular, once $C$ passes some threshold, its effect on the quality of the solution is minimum, and bigger values of $C$ only cause slightly larger numerical errors. This is also demonstrated in Table \ref{table:penaltyC}. 
For this problem, the viscosity solution is not smooth, so we do not expect the full $(k+1)$-th order accuracy for this example. However, for different values of $C$ ranging from 0.125 to 1.0, the numerical errors listed in Table \ref{table:penaltyC} are all of second order. Actually, for all of the simulations performed in this paper, we find that $C=0.25$ to be a good choice of the penalty constant. Unless otherwise noted, for the remaining of the paper, we will use $C=0.25$.
%\rm If instead, we take some positive values for the penalty constant $C$, then the viscosity solution could be obtained by our scheme. However, with a minimal value such as $C=0.001$, we are not adding enough viscosity to correct the entropy violation. See fig \ref{fig linnonsmth}.\\
\begin{figure}
  \centering
  \begin{tabular}{cc}
\includegraphics[width=.45\textwidth]{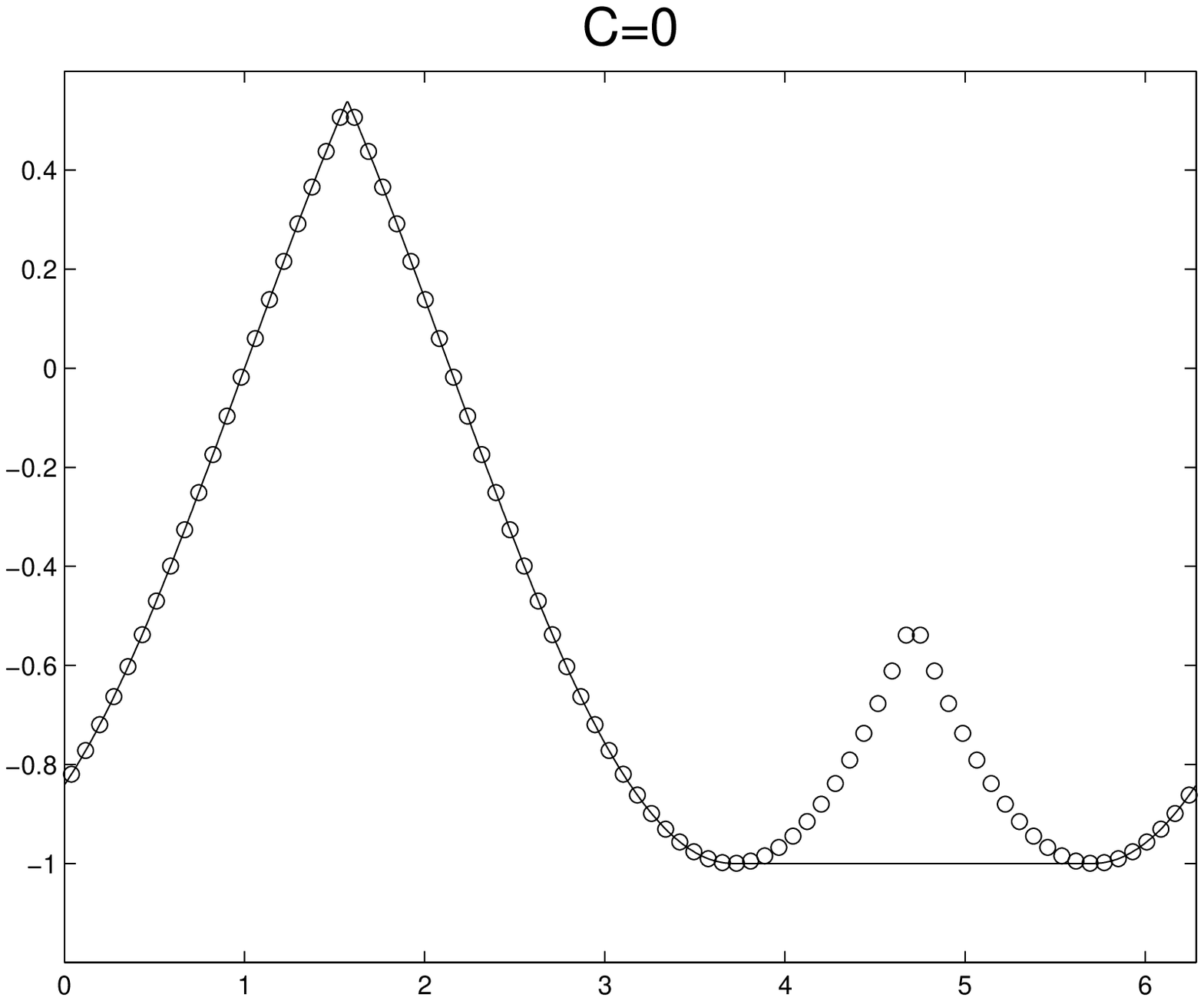}
\includegraphics[width=.45\textwidth]{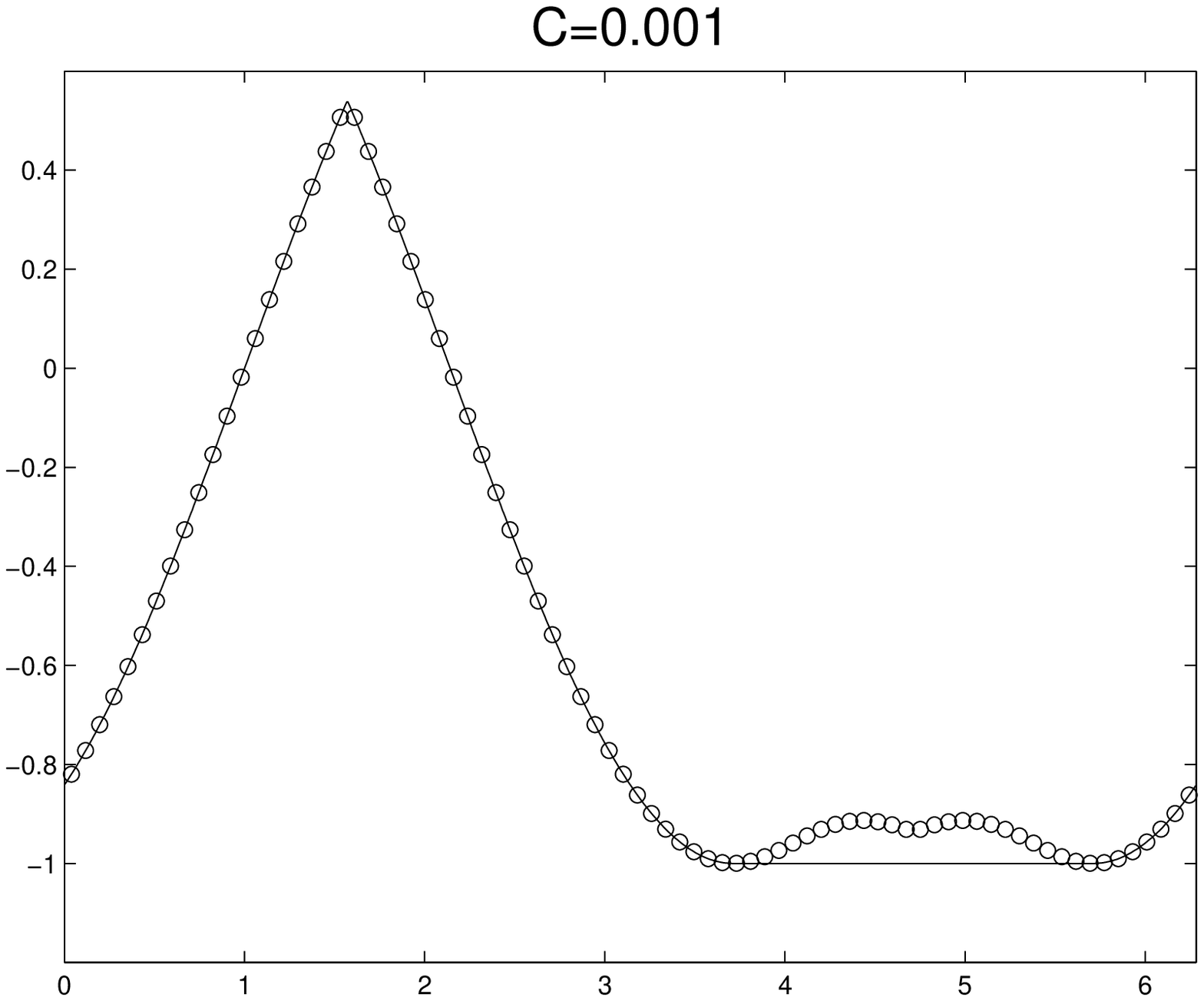}\\
\includegraphics[width=.45\textwidth]{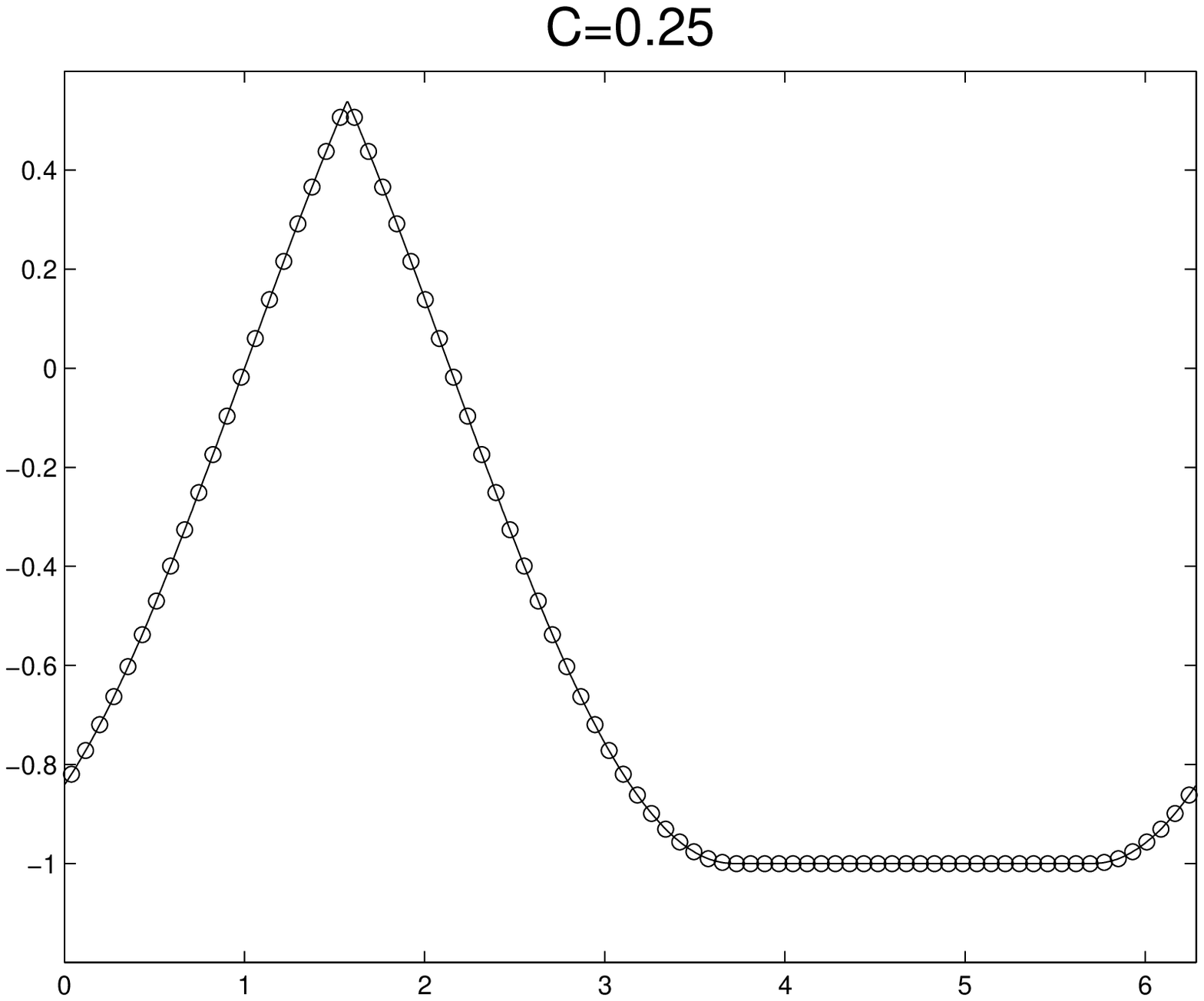}
\includegraphics[width=.45\textwidth]{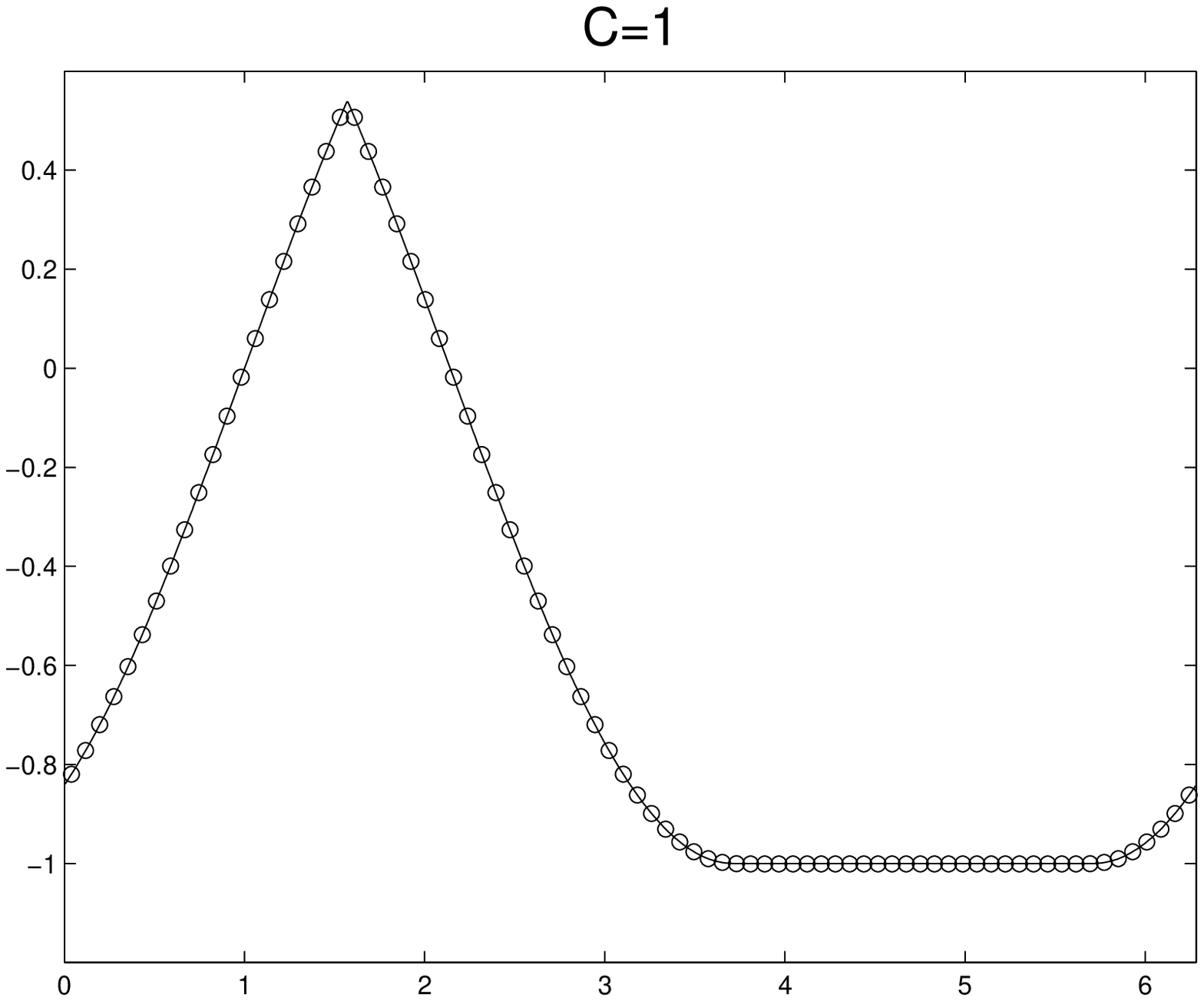}
\end{tabular}
  \caption{Example \ref{linnonsmth}. The numerical solution with various values of penalty constant $C$. $t=1$. $CFL=0.1$.   $P^2$ polynomials. $N=80$. Solid line: the exact solution; circles: the numerical solution. }
 \label{fig linnonsmth} 
\end{figure}

% \rm We could obtain the optimal order of accuracy using $P^1$ polynomials as shown in Table \ref{table:linnonsmth}.  However, we could not get the optimal order $k+1$ if we use $P^k$ polynomials for $k \ge 2$. The numerical errors are affected by the penalty term which introduces numerical errors of order $2$. The larger constant $C$ helps fix the entropy violation, however introduces more errors here.  The choice of proper constants $C$ will not affect the orders, which could be observed in Table \ref{table:penaltyC}.

%% Table 3.2
%\begin{table}[ht]
%\caption{Errors and numerical orders of accuracy for Example \ref{linnonsmth} when using $P^1$ polynomials and Runge-Kutta 
%third order time discretization on a uniform mesh of $N$ cells. Penalty constant $C=1.0$. Final time $t=1$. $CFL=0.1$.}
%\vspace{2 mm}
%\centering
%\begin{tabular}{c c c c c c c}
%\hline
%N   & $L^1$ error & order & $L^2$ error & order & $L^\infty$ error & order  \\
%\hline
%40 &  2.11E-03 &   & 2.97E-03 &    & 5.25E-03  &   \\
%80 &  5.22E-04 & 2.01 & 7.42E-04 & 2.00 & 1.32E-03 &1.99   \\
%160 & 1.30E-04 & 2.01 & 1.86E-04 & 1.99 & 3.50E-04 & 1.91  \\
%320 & 3.25E-05 & 2.00 &  4.68E-05 &  1.99 &1.09E-04 &  1.69  \\
%640 & 8.15E-06 & 2.00 & 1.18E-05 & 1.99 &3.26E-05  &  1.74 \\
%\hline
%\end{tabular}
%\label{table:linnonsmth}.
%\end{table}

% Table 3.3
\begin{table}[ht]
\caption{Errors and numerical orders of accuracy for Example \ref{linnonsmth} when using $P^2$ polynomials and third order Runge-Kutta 
 time discretization on a uniform mesh of $N$ cells.    Final time $t=1$. $CFL=0.1$.}
\vspace{2 mm}
\centering
\begin{tabular}{c c c c c c c}
\hline
N   & $L^1$ error & order & $L^2$ error & order & $L^\infty$ error & order  \\
\hline
 &  &      &   $C=1.0$ &   &  &   \\
 \hline
40 &  1.05E-03 &   & 1.85E-03 &   & 3.49E-03 &   \\
80 &  2.71E-04 & 1.95 & 4.78E-04 & 1.95 & 8.73E-04 &2.00   \\
160 & 6.89E-05 & 1.98 & 1.21E-04 & 1.98 & 2.18E-04 & 2.00  \\
320 & 1.73E-05 & 1.99 &  3.06E-05 &  1.98 &5.46E-05 &  2.00  \\
640 & 4.34E-06 & 2.00 & 7.67E-06 & 2.00 &1.37E-05  &  2.00  \\
\hline
 &  &      & $C=0.5$ &   &  &   \\
 \hline
40 &  9.92E-04 &   &  1.74E-03 &   & 3.28E-03&   \\
80 & 2.56E-04 & 1.95  & 4.50E-04 & 1.95 &8.22E-04 &2.00   \\
160 & 6.49E-05 & 1.98 & 1.14E-04 & 1.98 & 2.06E-04 & 2.00  \\
320 & 1.63E-05 & 1.99 & 2.88E-05 & 1.98  &5.14E-05 &  2.00  \\
640 & 4.09E-06 & 2.00 &7.22E-06 & 2.00 &1.29E-05  &  2.00  \\
\hline
 &  &      &  $C=0.25$ &   &  &   \\
\hline
40 &  8.74E-04 &   & 1.53E-03&   & 2.87E-03 &   \\
80 & 2.25E-04 & 1.96 & 3.95E-04 & 1.95 & 7.19E-04 &2.00   \\
160 & 5.69E-05 & 1.98 &1.00E-04 & 1.98 & 1.80E-04 & 2.00  \\
320 & 1.43E-05 &  1.99 &  2.52E-05 &  1.99 & 4.50E-05 &  2.00  \\
640 & 3.58E-06 & 2.00 & 6.32E-06 & 1.99 &1.13E-05  &  2.00  \\
\hline
 &  &      & $C=0.125$ &   &  &   \\
\hline
40 &  6.38E-04 &   & 1.10E-03 &   & 2.05E-03&   \\
80 & 1.62E-04 & 1.97 &  2.84E-04 & 1.96 & 5.14E-04 &2.00   \\
160 & 4.09E-05 & 1.98 & 7.18E-05 & 1.98 &  1.29E-04& 2.00  \\
320 & 1.03E-05 &  1.99 &  1.81E-05 &  1.99 & 3.23E-05 &  2.00  \\
640 & 2.57E-06 & 2.00 & 4.53E-06 & 2.00 & 8.57E-06  &  1.91 \\
\hline
\end{tabular}
\label{table:penaltyC}.
\end{table}

\end{exa}

%\subsection{Nonlinear smooth problems}

\begin{exa}\rm One-dimensional Burgers' equation with smooth initial condition
\label{nonlinsmth}
\begin{flalign}
\ \ \ \ \ & \begin{cases}
    \varphi_t+\frac{1}{2}\varphi_{x}^2=0  \\
    \varphi(x,0)=\sin(x)\\
    \varphi(0,t)=\varphi(2\pi,t)
  \end{cases}&
\end{flalign}
At $t=0.5$, the solution is still smooth, and the optimal $(k+1)$-th accuracy is obtained for $P^k$ polynomials with both uniform and random meshes, see Tables \ref{table:nonlinsmth} and   \ref{table:randommesh}. At $t=1$, there will be a shock in $\varphi_x$, and our scheme could capture the kink sharply as shown in Figure \ref{fig:nonlinsmth}.
 
% Table 3.4
\begin{table}[ht]
\caption{ Errors and numerical orders of accuracy for Example \ref{nonlinsmth} when using $P^2$ polynomials and third order Runge-Kutta 
 time discretization on a uniform mesh of $N$ cells. Penalty constant $C=0.25$. Final time $t=0.5$. $CFL=0.1$.}
\vspace{2 mm}
\centering
\begin{tabular}{c c c c c c c}
\hline
N   & $L^1$ error & order & $L^2$ error & order & $L^\infty$ error & order  \\
\hline
 &  & &$P^1$     &   & &   \\
\hline
40 &  8.45E-04 &   & 1.23E-03 &   & 5.04E-03&   \\
80 &  2.02E-04& 2.07 & 2.99E-04&2.04&  1.27E-03 & 1.99\\
160 &  4.93E-05 &  2.03 &7.42E-05 & 2.01 &  3.42E-04& 1.89 \\
320 & 1.22E-05 &2.01 &  1.86E-05 &  2.00 &9.08E-05 &  1.91\\
640 & 3.04E-06 &2.01 & 4.66E-06 &2.00 &2.36E-05&  1.94 \\
\hline
 &  &  & $P^2$    &   & &   \\
\hline
40 & 1.27E-05 &   & 2.33E-05&   & 1.28E-04&   \\
80 & 1.53E-06 & 3.05 &  2.93E-06 & 2.99 & 2.10E-05 & 2.61 \\
160 & 1.91E-07 & 3.00 & 3.73E-07 & 2.98 &  2.52E-06 &  3.06  \\
320 &  2.39E-08 &  3.00 &  4.74E-08 &  2.98 & 3.56E-07&  2.82  \\
640 &3.63E-09 & 2.72  &6.23E-09 &2.93 & 4.82E-08 &  2.88 \\
\hline
\end{tabular}
\label{table:nonlinsmth}
\end{table}

% random mesh
\begin{table}[ht]
\caption{ Errors and numerical orders of accuracy for Example  \ref{nonlinsmth}  when using $P^1$ and $P^2$ polynomials and third order Runge-Kutta 
time discretization on a random mesh with $40\%$ perturbation of $N$ cells. Penalty constant $C=0.25$. Final time $t=0.5$. $CFL=0.1$.}
\vspace{2 mm}
\centering
\begin{tabular}{c c c c c c c}
\hline
N   & $L^1$ error & order & $L^2$ error & order & $L^\infty$ error & order  \\
\hline
 &  & &$P^1$     &   & &   \\
\hline
40 &  1.23E-03 &   &  1.91E-03 &   & 1.01E-02&   \\
80 &  2.70E-04 & 2.19 &4.25E-04&2.17 &  2.59E-03 &1.96\\
160 &  6.70E-05 &  2.01 & 1.05E-04 & 2.01 &  6.22E-04& 2.06 \\
320 & 1.62E-05 & 2.05 &  2.67E-05 &  1.97 &2.03E-04 &  1.61\\
640 &  3.97E-06 &2.03 & 6.69E-06 &2.00 &6.52E-05 & 1.64 \\
\hline
 &  &  & $P^2$    &   & &   \\
\hline

40 &  2.27E-05 &   &  4.52E-05 &   & 2.96E-04&   \\
80 &  2.54E-06 &  3.16 &5.84E-06 & 2.95 & 5.25E-05 &2.50 \\
160 &  3.19E-07 & 3.00 & 6.87E-07 & 3.09  & 5.82E-06& 3.17\\
320 &  4.00E-08 &  3.00 &  9.34E-08 &  2.88 &8.96E-07 & 2.70  \\
640 & 5.38E-09 &  2.89  &  1.16E-08 &3.01 & 1.32E-07&2.77 \\
\hline
\end{tabular}
\label{table:randommesh}
\end{table}

% graph
\begin{figure}[h!]  
\centering
    \includegraphics[width=0.6\textwidth]{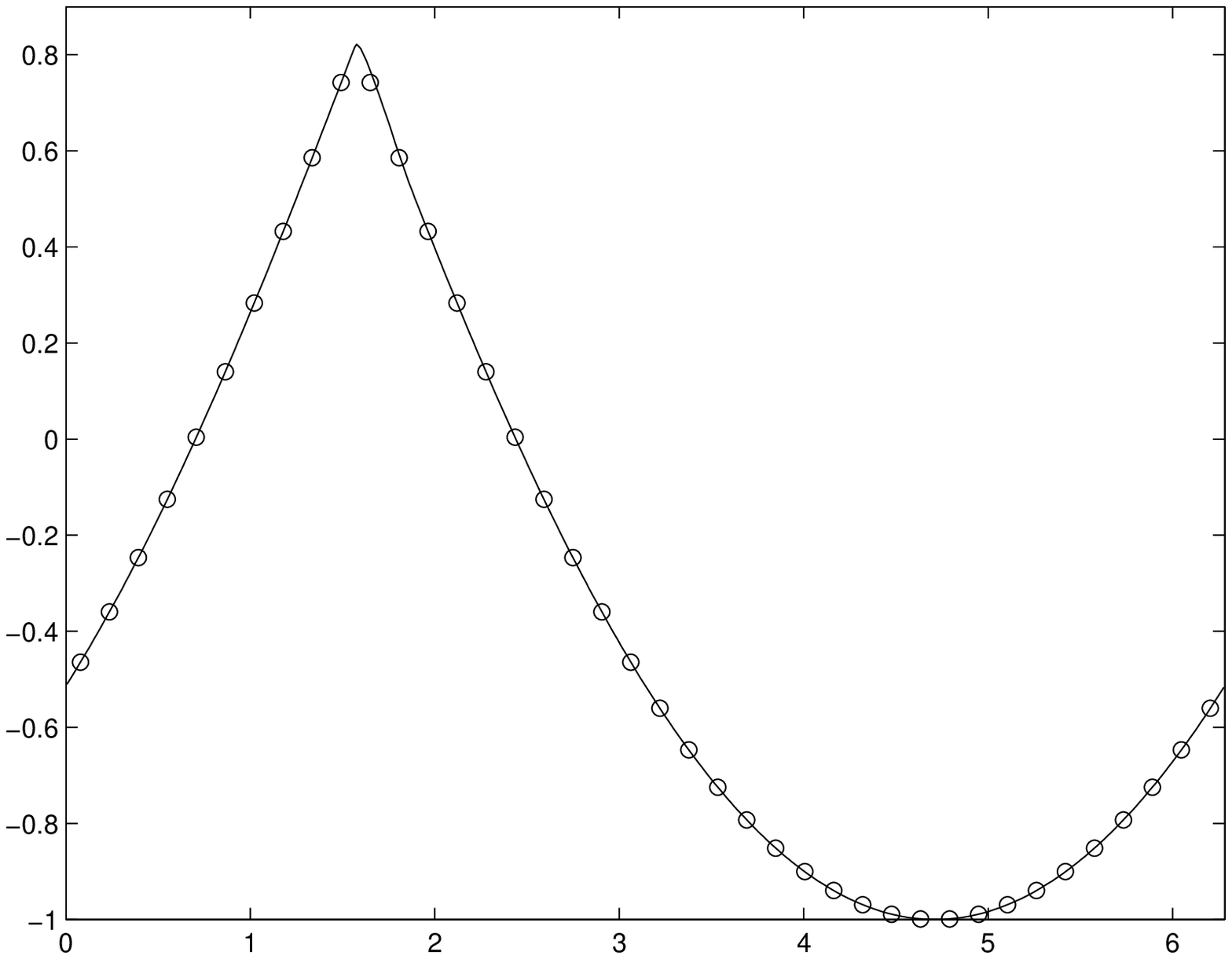}
\caption{Example \ref{nonlinsmth}. $t=1.5$. $CFL=0.1$.  $P^2$ polynomials. $N=40$. Penalty constant $C=0.25$. Solid line: the exact solution; circles: the numerical solution.} 
\label{fig:nonlinsmth}
\end{figure}
\end{exa}

\begin{exa}\rm One-dimensional Burgers' equation with a nonsmooth initial condition  \cite{Cheng_07_JCP_HJ}
\label{bgnonsmthIC}
\begin{flalign}
\ \ \ \ \ & \begin{cases}
    \varphi_t+\frac{1}{2}\varphi_{x}^2=0  \\
    \varphi(x,0)= 
    \begin{cases}
    \pi-x & \textrm{if} \quad 0 \le x \le \pi\\
    x-\pi & \textrm{if} \quad 0\le x \le 2 \pi\\
    \end{cases}\\   
     \varphi(0,t)=\varphi(2\pi,t) 
  \end{cases}&
\end{flalign}
The exact solution should have a rarefaction wave forming in its derivative, so the initial sharp corner at $x=\pi$ should be smeared out at later times. Since the entropy condition is violated by the Roe type scheme,  the entropy fix is necessary for convergence. Figure \ref{fig:bgnonsmthIC} shows the comparison of our schemes with  various values of penalty constant $C$ for this nonlinear problem. Again, we could see that $C=0.25$ is a good choice for this example.

\begin{figure}
  \centering
  \begin{tabular}{cc}
\includegraphics[width=.45\textwidth]{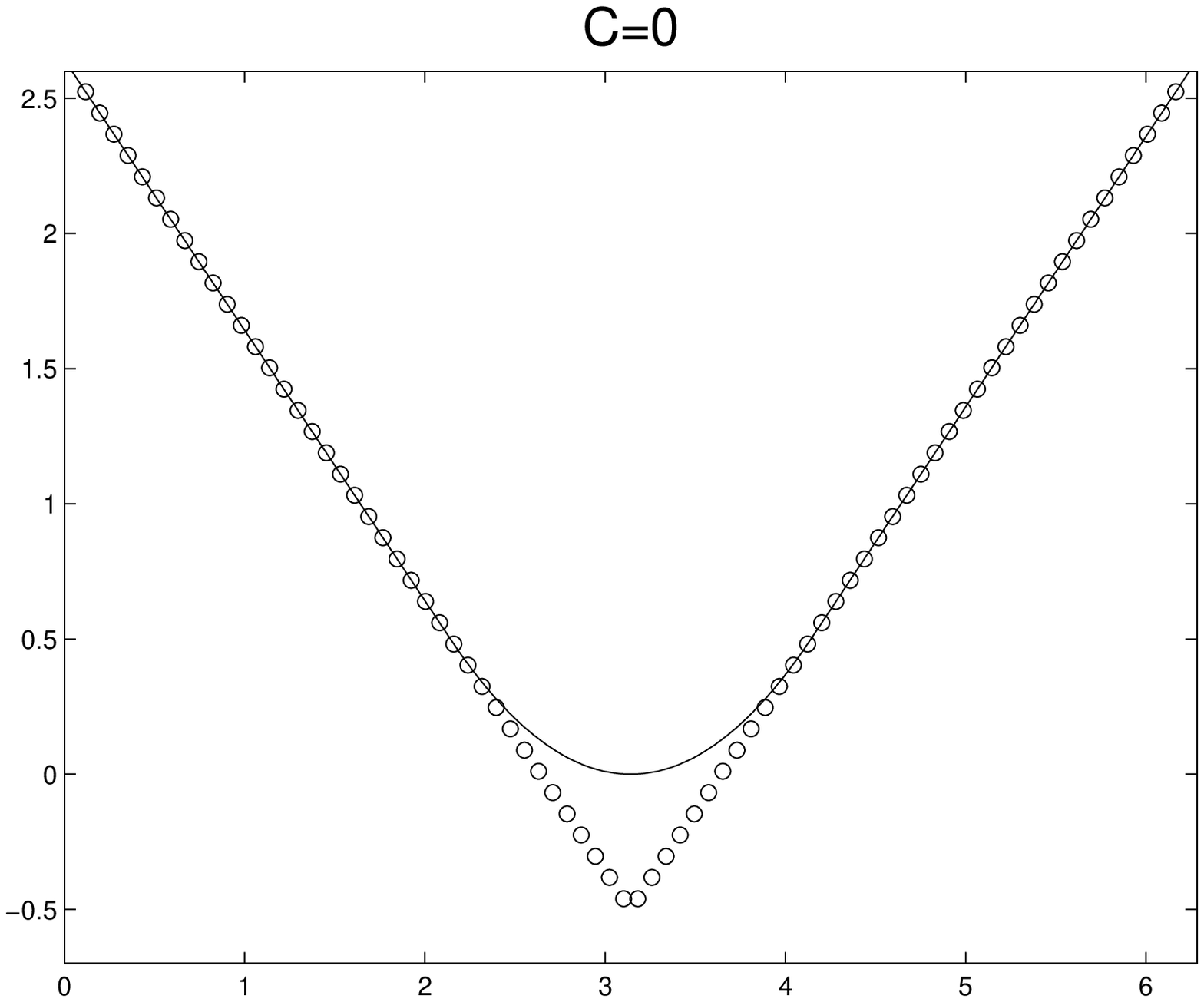}
\includegraphics[width=.45\textwidth]{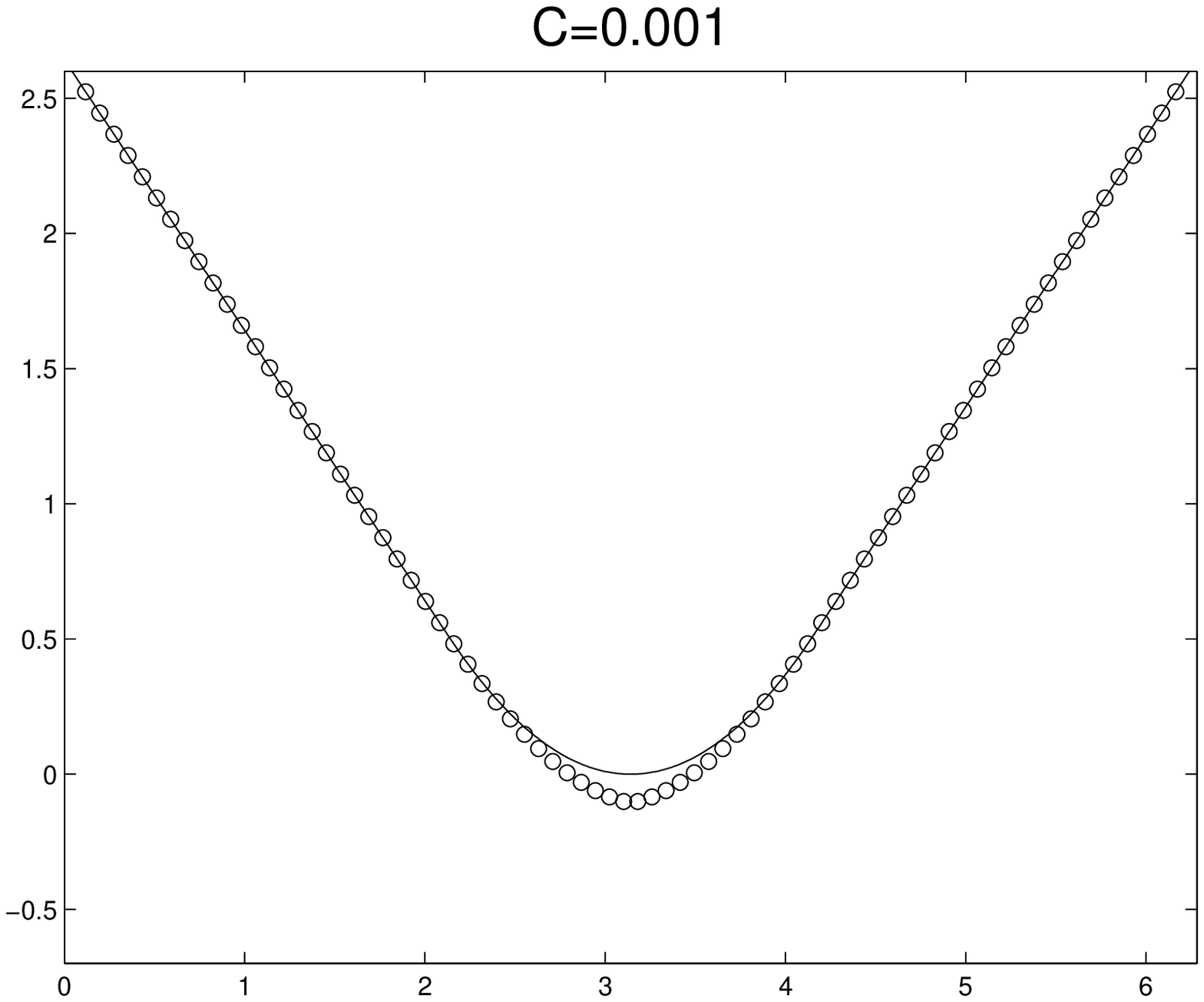}\\
\includegraphics[width=.45\textwidth]{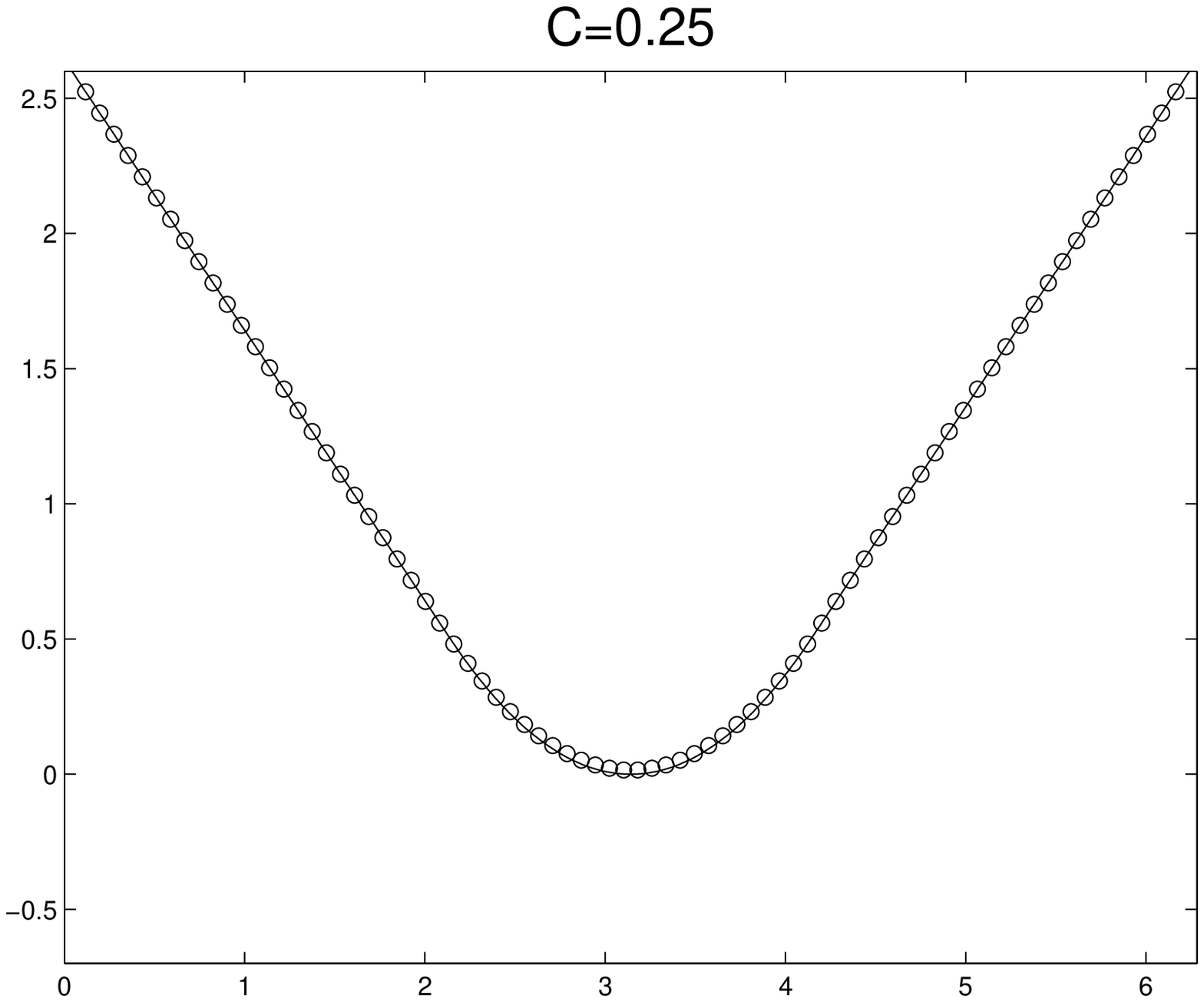}
\includegraphics[width=.45\textwidth]{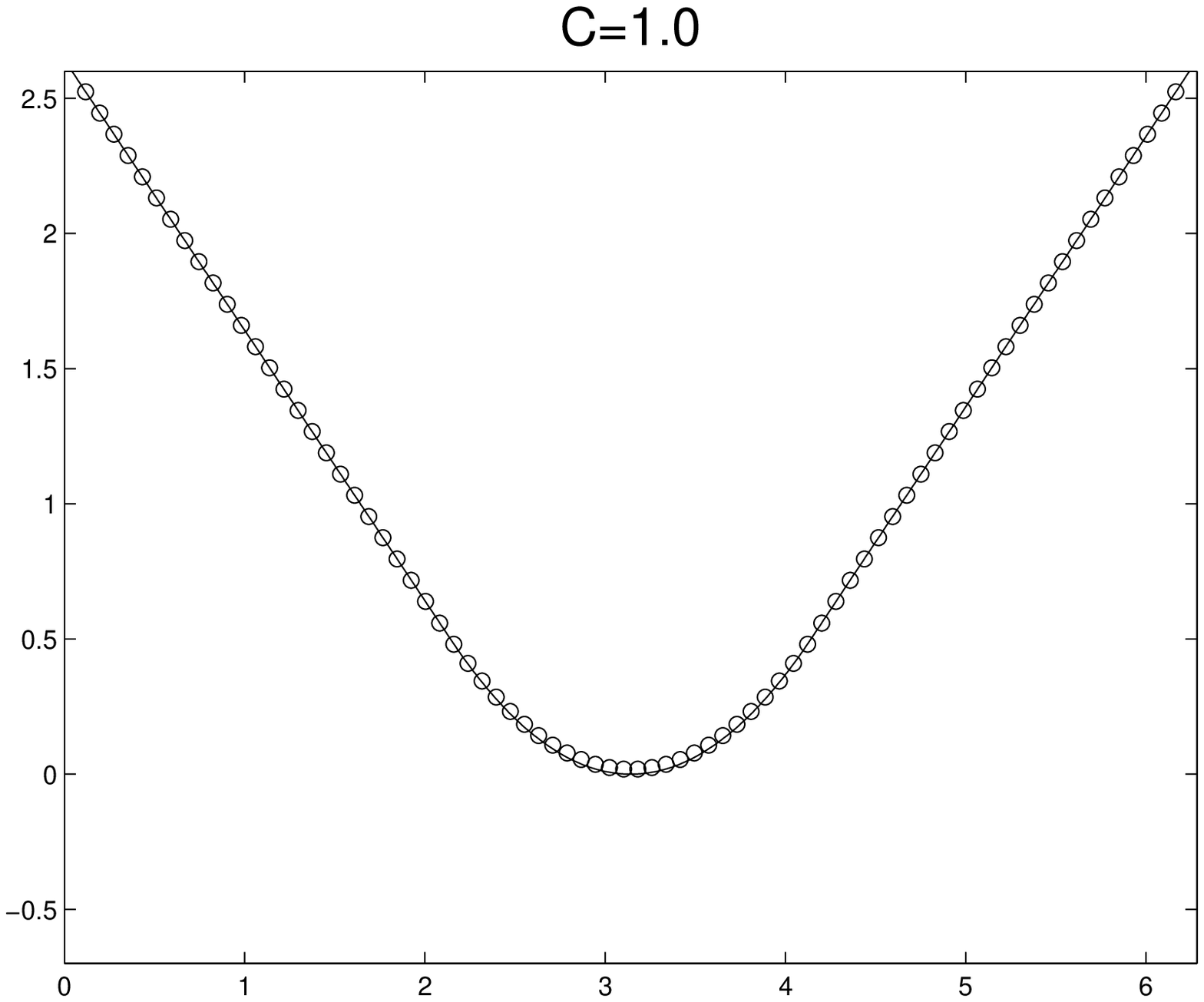}
\end{tabular}
  \caption{Example \ref{bgnonsmthIC}. The numerical solution with various values of penalty constant $C$.  $t=1$. $CFL=0.1$.   $P^2$ polynomials. $N=80$.   Solid line: the exact solution; circles: the numerical solution. }
 \label{fig:bgnonsmthIC} 
\end{figure}

%\begin{figure}
%
%\centering
%\begin{minipage}{.4\textwidth}
%  \centering
%  \includegraphics[width=\linewidth]{nocorrection}
%\end{minipage}
%\begin{minipage}{.4\textwidth}
% \centering
%  \includegraphics[width=\linewidth]{correction}
%\end{minipage}
%\caption{Example \ref{bgnonsmthIC}. $t=1$. $CFL=0.05$. using $P^2$ polynomials. $N=40$. Solid line: the exact solution; circles: the numerical solution. Left:  without entropy fix $C=0$; Right: with entropy fix $C=0.25$. }
%\label{fig:bgnonsmthIC}
%\end{figure}

\end{exa}

%\subsection{Nonlinear nonsmooth problems}

\begin{exa}\rm One-dimensional Eikonal equation \label{Eikonal}
\begin{flalign}
 \ \ \ \ \  &\begin{cases}
    \varphi_t+|\varphi_{x}|=0 \\
    \varphi(x,0)=\sin(x) \\
    \varphi(0,t)=\varphi(2\pi,t) 
  \end{cases}&
\end{flalign}
The exact solution is the same as the exact solution of Example \ref{linnonsmth}. 
%The Hamiltonian is non smooth and convex. The entropy correction is automatically turned on when needed in the entropy violating cells. 
Our scheme could capture the viscosity solution of this nonsmooth Hamiltonian. The numerical errors and orders of accuracy using $P^2$ polynomials are listed in Table \ref{table:Eikonal}. Since the solution is not smooth, we do not expect the optimal $(k+1)$-th order accuracy for $P^k$ polynomials.
% Table 3.5
\begin{table}[ht]
\caption{ Errors and numerical orders of accuracy for Example \ref{Eikonal} when using $P^2$ polynomials and third order  Runge-Kutta 
time discretization on a uniform mesh of $N$ cells. Penalty constant $C=0.25$. Final time $t=1$. $CFL=0.1$.}
\vspace{2 mm}
\centering
\begin{tabular}{c c c c c c c}
\hline
N   & $L^1$ error & order & $L^2$ error & order & $L^\infty$ error & order  \\
\hline
40 &  6.24E-04 &   & 1.09E-03 &   &2.13E-03&   \\
80 &  1.69E-04 & 1.88 &  2.98E-04 & 1.87 & 5.54E-04 &1.94 \\
160 & 4.35E-05 & 1.96 & 7.67E-05 & 1.96 &  1.40E-04 &  1.98 \\
320 & 1.10E-05 &  1.99 &   1.94E-05 &  1.99 & 3.51E-05 &  2.00  \\
640 & 2.75E-06 & 2.00  & 4.88E-06 &1.99 & 8.77E-06 &  2.00\\
\hline
\end{tabular}
\label{table:Eikonal}
\end{table}
\end{exa}

%\subsection{Nonconvex problems}
\begin{exa}
\label{ex:nonconvsm} \rm One-dimensional equation with a nonconvex Hamiltonian
\begin{flalign}
\ \ \ \ \ &\begin{cases}
\varphi_t - \cos (\varphi_x + 1) = 0 \\
\varphi(x,0) = - \cos(\pi x) \\
    \varphi(-1,t)=\varphi(1,t)
\end{cases}&
\end{flalign}

\rm This example involves a nonconvex Hamiltonian with  smooth initial data. At $t = 0.5/\pi^2$, the exact solution is still smooth, and numerical results are presented in Table \ref{table:nonconvsm}, demonstrating the optimal order of accuracy of the scheme. By the time $t = 1.5/\pi^2$, nonsmooth features would develop in  $\varphi$, which are reliably captured in Figure \ref{fig:nonconvsm}.

% Table 3.6
\begin{table}[ht]
\caption{ Errors and numerical orders of accuracy for Example \ref{ex:nonconvsm}  when using $P^2$ polynomials and third order Runge-Kutta 
time discretization on a uniform mesh of $N$ cells. Penalty constants $C=0.25$. Final time $t=0.5/\pi^2$. $CFL=0.1$.}
 
\vspace{2 mm}
\centering
\begin{tabular}{c c c c c c c}
\hline
N   & $L^1$ error & order & $L^2$ error & order & $L^\infty$ error & order  \\
 \hline
40 &   1.46E-05 &   & 2.16E-05 &   &  9.89E-05 &   \\
80 &  1.79E-06 &  3.02 &2.87E-06& 2.91 & 1.59E-05 &2.64   \\
160 &2.22E-07 &  3.01 & 3.73E-07 &  2.95 & 2.39E-06 &2.74 \\
320 & 2.76E-08 & 3.01 &  4.79E-08 &   2.96 & 3.39E-07& 2.82 \\
640 & 3.51E-09 &2.98 & 6.13E-09 &  2.97&4.53E-08  &  2.90 \\
\hline
\end{tabular}
\label{table:nonconvsm}
\end{table}

% graph
\begin{figure}[h!]
\centering
    \includegraphics[width=0.6\textwidth]{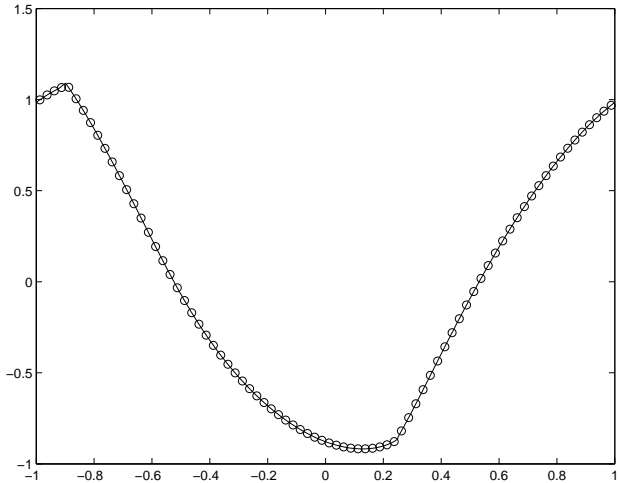}
 \label{fig:nonconve}
 \caption{Example \ref{ex:nonconvsm}. $t=1.5/\pi^2$. $CFL=0.1$.   $P^2$ polynomials. $N=80$. Penalty constant $C=0.25$. Solid line: the exact solution; circles: the numerical solution. }
 \label{fig:nonconvsm}
\end{figure}
\end{exa}

\begin{exa}
\label{ex:nonconvr} 
\rm  One-dimensional Riemann problem with a nonconvex Hamiltonian
\begin{flalign}
\ \ \ \ \ &\begin{cases}
\varphi_t + \frac{1}{4}(\varphi_x^2 - 1)(\varphi_x^2 - 4) = 0 \\
\varphi(x,0) = - 2 |x|
\end{cases}&
\end{flalign}\\
 For this problem, the initial condition has a singularity at $x=0$.
 Similar to \cite{li2010central,yan2011local}, a nonlinear limiter is needed in order to capture the viscosity solution.  We use the standard minmod limiter \cite{Cockburn_2001_RK_DG}. This example and Example \ref{2dReimann} are the only examples needing nonlinear limiting in this paper. 
 
 The numerical solutions with and without the limiter  are listed in Figure \ref{riemann2} for odd and  even values of $N$. Those different behaviors are due to the fact that the singular point $x=0$ would be exactly located at the cell interface for even $N$ but not odd $N$ at $t=0$. We note that the method with limiter can correctly capture the viscosity solution for both even and odd $N$.
 The numerical errors and orders of accuracy using $P^2$ polynomials with limiters are listed in Table \ref{nonconvex}. We could see that both methods converge, while the odd $N$ giving slightly smaller errors. However, similar to \cite{li2010central}, the method is only first order accurate for this nonsmooth problem.
 
% Given the initial condition is not smooth at $x=0$, the numerical results would depend on whether $N$ is even or odd. Since the information our scheme demands lies on the edges of the cells, odd numbers of cells would avoid the situation that the discontinuity point lies on the cell boundaries. Without the limiter, the numerical scheme would not revolve as the time goes when $N$ is even. Take an odd number $N=81$ for example, a comparison of the numerical results with and without the limiter can be seen from fig $3.5$. The minmod limiter is applied where the exact solution is not linear. The computation is carried out with the outflow numerical boundary conditions. 
% Odd meshes give better approximations actually, but need smaller CFL numbers. Smaller errors if N=161, 
%L-1 error =    6.41E-04, L-2 error =1.61E-03, L-infinity error = 9.88E-03

%\rm The $P^1$ and $P^2$ approximations will have comparable spatial orders for this example due to the usage of the nonlinear limiter. The numerical errors and orders of accuracy using $P^2$ polynomials are listed in Table \ref{nonconvex}. It could clearly be seen that with finer meshes, our scheme could better capture the exact solution, which is demonstrated by fig $3.6$. 
\begin{figure}
  \centering
  \begin{tabular}{cc}
\includegraphics[width=.45\textwidth]{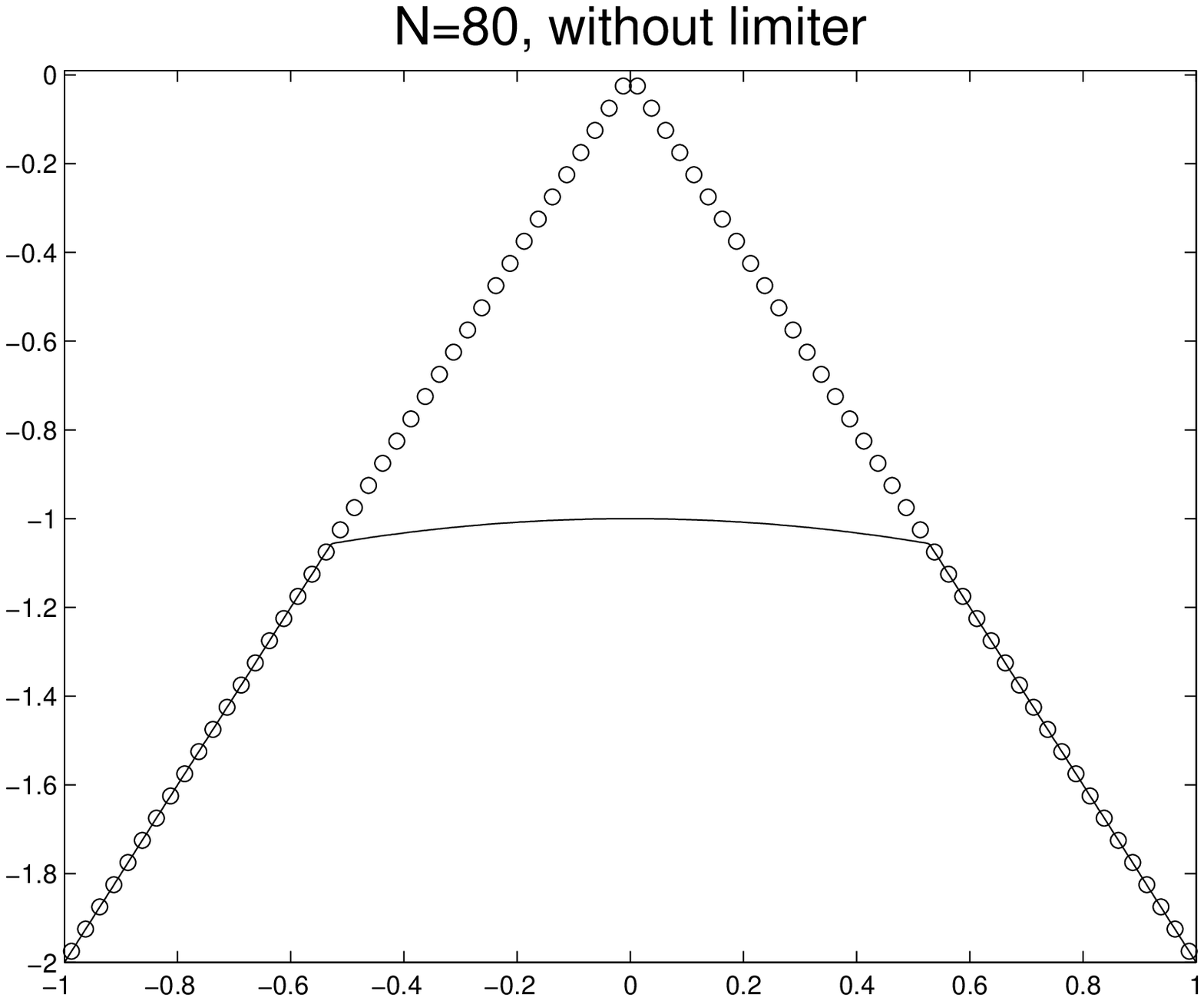}
\includegraphics[width=.45\textwidth]{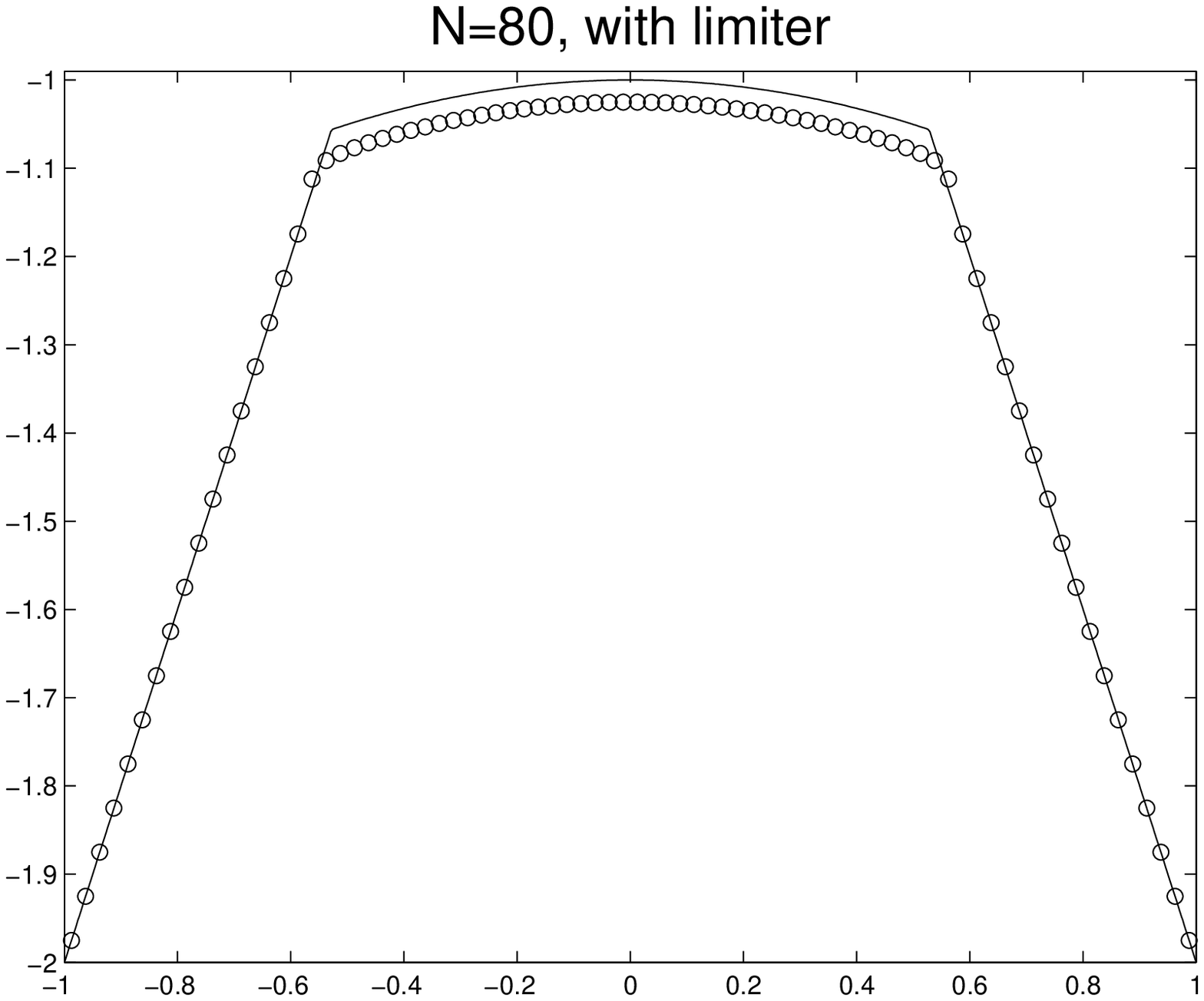}\\
\includegraphics[width=.45\textwidth]{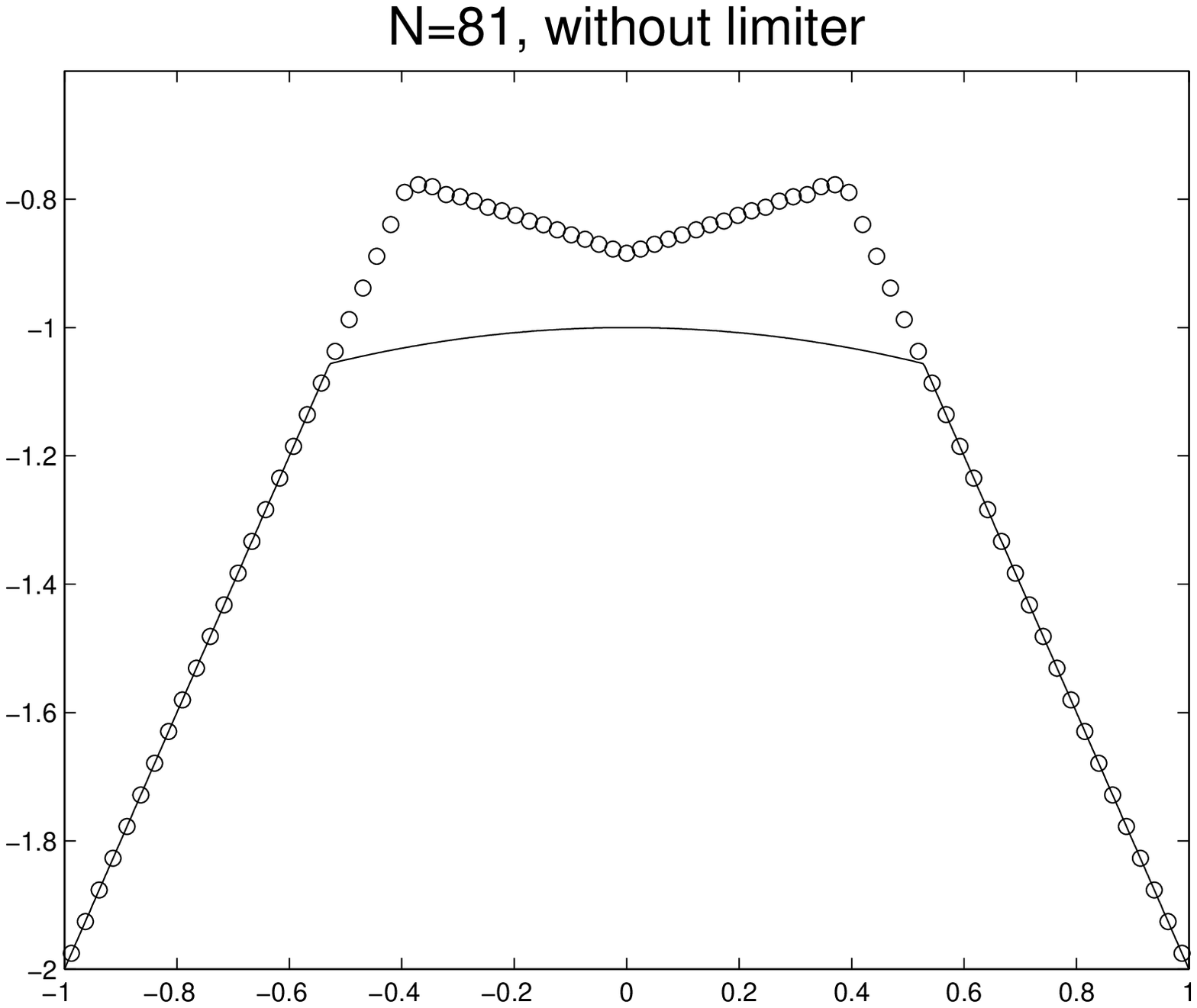}
\includegraphics[width=.45\textwidth]{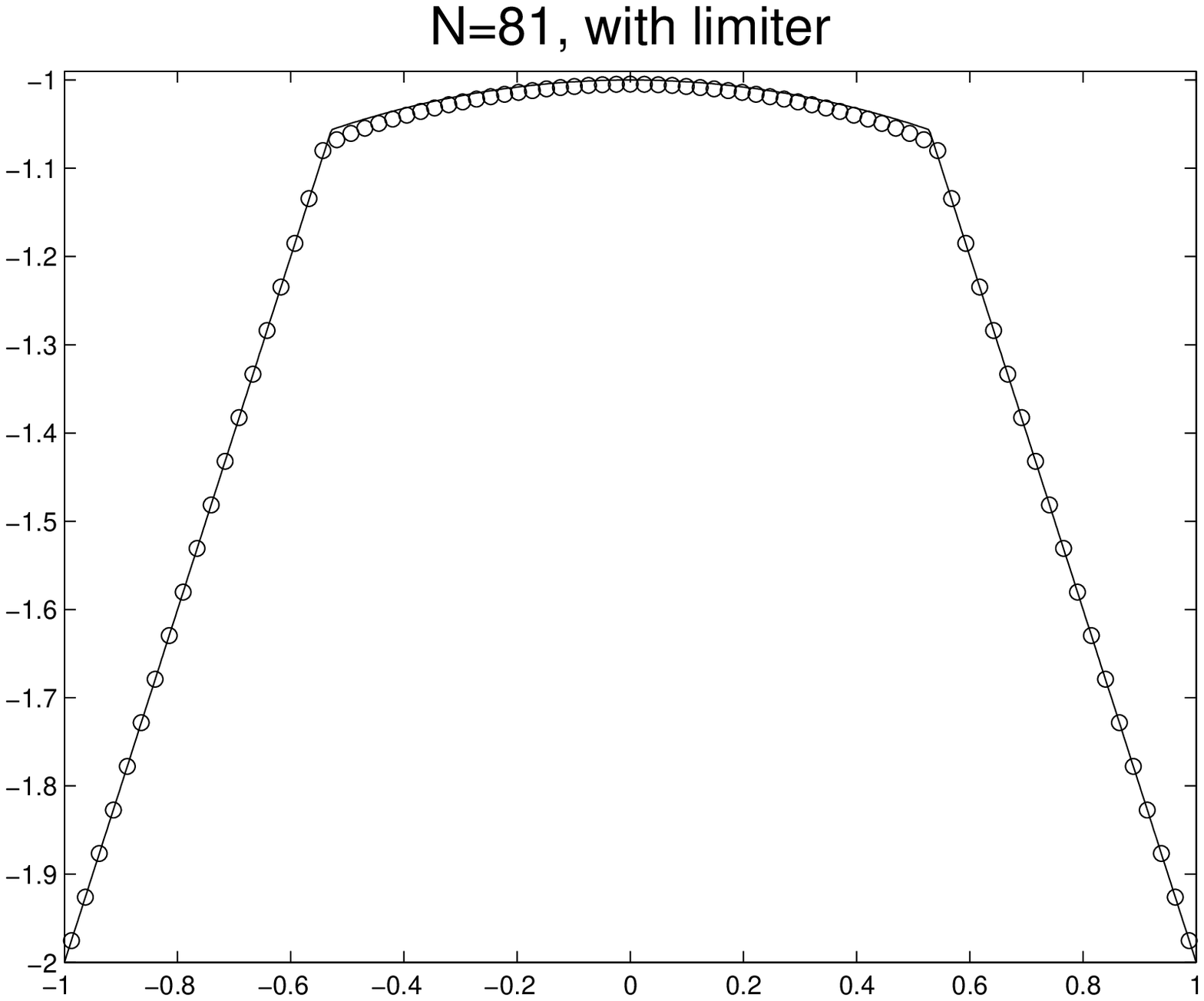}
\end{tabular}
  \caption{Example \ref{ex:nonconvr}. Comparison of the numerical solution with and without the limiter. $t=1$.   $P^2$ polynomials. $CFL=0.05$. Penalty constant $C=0.25$.  Left: without limiter; right: with limiter. Solid line: the exact solution; circles: the numerical solution. }
 \label{riemann2} 
\end{figure}

% Change
% Table 3.7
\begin{table}[ht]
\caption{ Errors and numerical orders of accuracy for Example \ref{ex:nonconvr}  when using $P^2$ polynomials and third order  Runge-Kutta 
time discretization on a uniform mesh of $N$ cells.  $CFL=0.05$. Penalty constant $C=0.25$. Final time $t=1$. A minmod limiter is used. }
\vspace{2 mm}
\centering
\begin{tabular}{c c c c c c c}
\hline
N   & $L^1$ error & order & $L^2$ error & order & $L^\infty$ error & order  \\
\hline
 &    &  &Even N  &   & &   \\
\hline
40 &   9.49E-03 &   & 2.21E-02 &   & 5.96E-02&   \\
80 & 4.64E-03 & 1.03 &    1.10E-02 &1.00 &3.17E-02 &0.91 \\
160 & 2.28E-03 & 1.03 & 5.48E-03 & 1.00 &  1.64E-02 &  0.95 \\
320 & 1.12E-03 &  1.02 &   2.73E-03 &  1.01& 8.40E-03 &  0.97  \\
640 &   5.60E-04 & 1.00  & 1.36E-03 &1.00 & 4.27E-03 & 0.98\\
\hline
 &    & & Odd N  &    & &   \\
 \hline
 41 &   2.81E-03 &   & 6.74E-03 &   & 2.94E-02&   \\
81 & 1.34E-03  & 1.09 &    3.35E-03 &1.03 & 2.38E-02 &0.31 \\
161 & 6.41E-04 & 1.07 &  1.61E-03 & 1.06 &   9.88E-03 & 1.28 \\
321 & 3.17E-04 &  1.02 &    7.99E-04 &  1.01& 4.36E-03 &  1.19  \\
641 &   1.56E-04 &  1.02  & 3.96E-04 &1.02 & 3.12E-03 & 0.49\\
\hline
\label{nonconvex}
\end{tabular}
\end{table}

%\begin{figure}
%\centering
%\begin{minipage}{.4\textwidth}
%  \centering
%  \includegraphics[width=\linewidth]{nolimP281}
%\end{minipage}
%\begin{minipage}{.4\textwidth}
% \centering
%  \includegraphics[width=\linewidth]{P281}
%\end{minipage}
%\caption{Example \ref{ex:nonconvr}  $t=1$. $CFL=0.05$. using $P^2$ polynomials. $N=81$. Left: without limiter; Right: with limiter. Solid line: the exact solution; circles: the numerical solution. }
%\label{riemann1}
%\end{figure}
%
%\begin{figure}
%\centering
%\begin{minipage}{.4\textwidth}
%  \centering
%  \includegraphics[width=\linewidth]{P280nolim}
%\end{minipage}
%\begin{minipage}{.4\textwidth}
% \centering
%  \includegraphics[width=\linewidth]{P2160}
%\end{minipage}
%\caption{Example \ref{ex:nonconvr}  $t=1$. $CFL=0.1$. using $P^2$ polynomials. Left: $N=80$; Right: $N=160$. Solid line: the exact solution; circles: the numerical solution. }
%\label{riemann2}
%\end{figure}

\end{exa}

%\clearpage
\subsection{Two-dimensional results}
In this subsection, we provide computational results for two-dimensional HJ equations on both Cartesian and unstructured meshes.

% Cheng 4.1.2
\begin{exa}\rm Two-dimensional linear problem with smooth variable coefficients
\label{rotating.hill}
\begin{equation}
\label{rotating.hill.eqn}
\varphi_t-y\varphi_x+x\varphi_y=0 .
\end{equation}
The computational domain is $[-1,1]^2$. The initial condition is given by
\begin{equation}
\label{rotating.hill.init}
\varphi_0(x,y) = \left\{
      \begin{array} {l l}
       0 & 0.3 \leq r \\
       0.3-r & 0.1 < r < 0.3 \\
       0.2   & r \leq 0.1	
      \end{array}
      \right.
\end{equation}
where $r=\sqrt{(x-0.4)^2+(y-0.4)^2}$. 
We   impose periodic boundary condition on the domain. 
This is a solid body rotation around the origin.
The exact solution can be expressed as
\begin{equation}
\label{rotating.hill.exact}
\varphi(x,y,t)=\varphi_0(x\cos(t)+y\sin(t),-x\sin(t)+y\cos(t)).
\end{equation}

For this problem, same as the argument in Example \ref{linsmth}, the choice of $C$ does not have an effect on the scheme. We list the numerical errors and orders in Table \ref{table hill.p2}. With this nonsmooth initial condition, we do not expect to obtain ($k+1$)-th order of accuracy. At $t=2 \pi$, i.e. one period of rotation, we take a snapshot at 
the line $y=x$ in Figure \ref{fig hill.p2}.  It can be clearly seen that a higher order 
scheme can yield better results for this nonsmooth initial condition.

%Table 3.8
\begin{table}[ht]
\caption {Errors and numerical orders of accuracy for 
Example  \ref{rotating.hill} when using $P^2$ polynomials and third order Runge-Kutta 
 time discretization on a uniform mesh of $N \times N$ cells. 
Final time $t=1$. $CFL=0.1$.}
\vspace{2 mm}
\centering
\begin{tabular}{c c c c c c c}
\hline
N   & $L^1$ error & order & $L^2$ error & order & $L^\infty$ error & order  \\
\hline
 10 & 1.21E-03 &  & 3.10E-03 &   & 2.21E-02 & \\
 20 & 4.13E-04 &   1.55   & 1.32E-03 &     1.23  & 1.14E-02 &0.95   \\
 40 & 1.38E-04 & 1.58 &5.51E-04 & 1.26 & 6.49E-03 &   0.81  \\
80 & 4.74E-05 & 1.54 &  2.36E-04 &  1.22 & 3.62E-03&   0.84  \\
160 &1.54E-05& 1.62&1.01E-04&1.23 & 2.07E-03&0.81\\
\hline
\end{tabular}
\label{table hill.p2}.
\end{table}

\begin{figure}
\centering
\begin{minipage}{.5\textwidth}
  \centering
  \includegraphics[width=\linewidth]{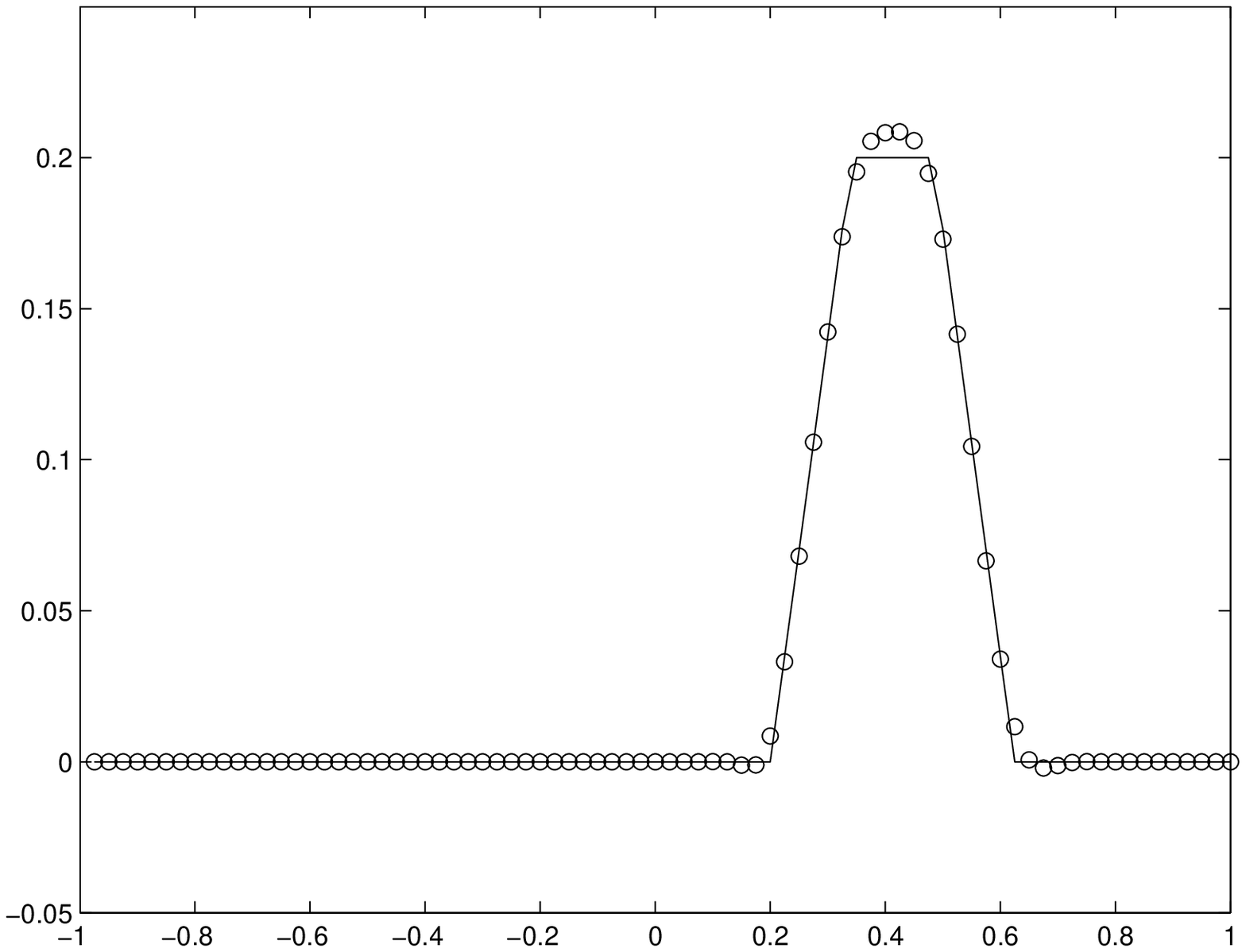}
\end{minipage}%
\begin{minipage}{.5\textwidth}
 \centering
  \includegraphics[width=\linewidth]{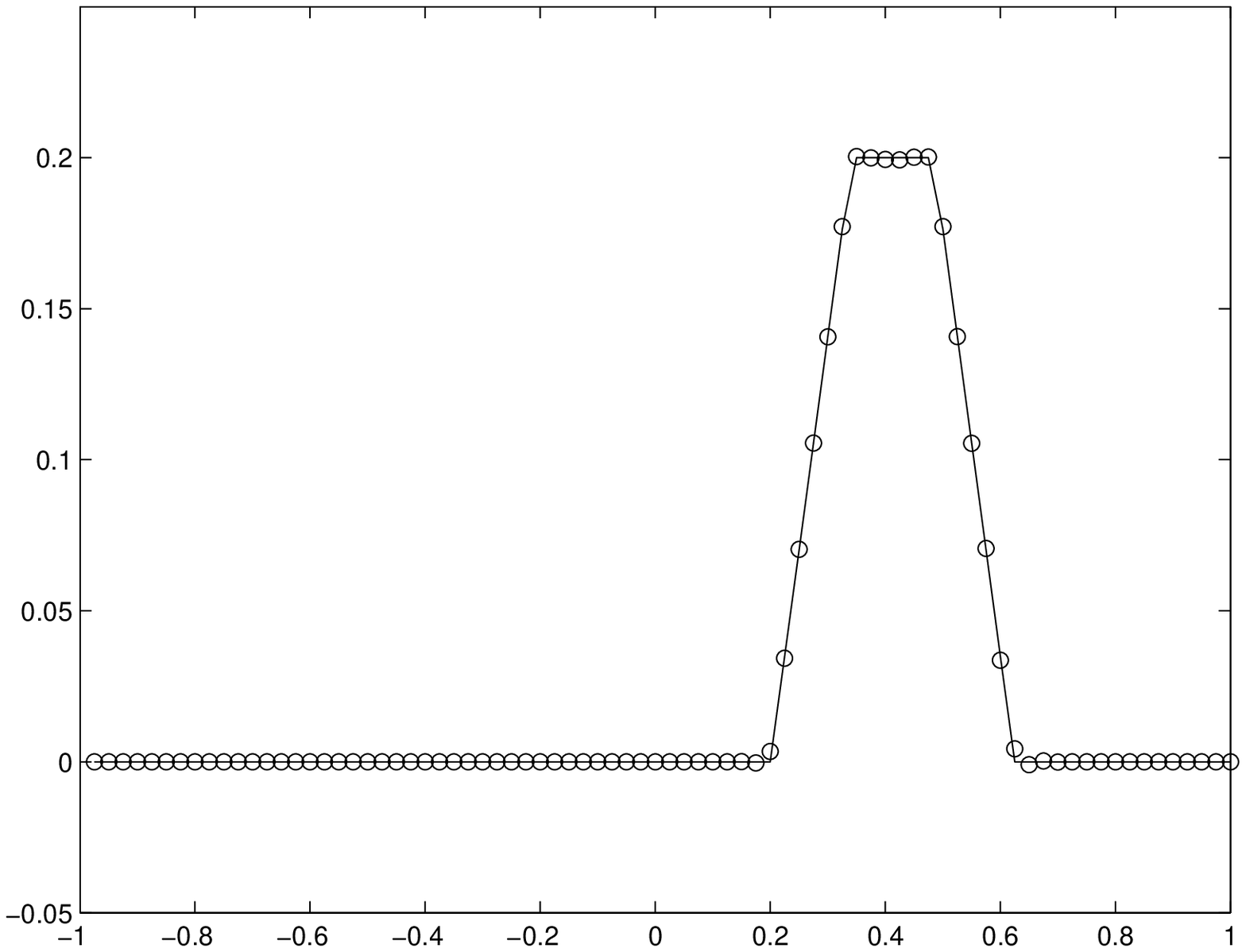}
\end{minipage}
\caption{Example \ref{rotating.hill}. $t=2\pi$. $CFL=0.1$.  $80\times 80$ uniform mesh. Left: $P^1$ polynomials; right: $P^2$ polynomials. One dimensional cut of $45^\circ$ with the $x$ axis. Solid line: the exact solution; circles: the numerical solution. }
\label{fig hill.p2}
\end{figure}
\end{exa}

\begin{exa}\rm We solve the same equation (\ref{rotating.hill.eqn}) as  
in Example \ref{rotating.hill}, but with a smooth initial condition as
\label{rotating.hill.smoothint}
\begin{equation}
\varphi_0(x,y) = \exp \left( -\frac{(x-0.4)^2+(y-0.4)^2}{2 \sigma^2} \right).
\end{equation}
The constant $\sigma=0.05$ is chosen such that at the domain boundary, $\varphi$ is 
very small, hence imposing the periodic boundary condition will lead
to small  errors.  We then could observe the optimal order 
of accuracy in Table \ref{table hillsmooth.p2}.

\begin{table}[ht]
\caption {Errors and numerical orders of accuracy for 
Example 3.9 when using $P^2$ polynomials and third order Runge-Kutta 
 time discretization on a uniform mesh of $N \times N$ cells. 
Final time $t=1$. $CFL=0.1$.}
\vspace{2 mm}
\centering
\begin{tabular}{c c c c c c c}
\hline
N   & $L^1$ error & order & $L^2$ error & order & $L^\infty$ error & order  \\
\hline
 20 & 1.42E-03 &   & 1.03E-02 &  &  2.79E-01 &  \\
 40 & 1.54E-04 & 3.20 &1.47E-03 &2.81 &5.25E-02 &  2.41  \\
80 & 1.10E-05 & 3.81 & 1.10E-04&  3.73 & 5.77E-03 & 3.19  \\
160 & 1.12E-06 & 3.30 & 1.15E-05&  3.26 &  8.96E-04 &   2.69 \\
\hline
\end{tabular}
\label{table hillsmooth.p2}
\end{table}
\end{exa}

%triangle code in /chengcode/triangles
\begin{exa}
\label{2dburgtr}
\rm  Two-dimensional Burgers' equation
\begin{equation}
\label{2d.burg} \left\{\begin{array} {l}
        \displaystyle    \varphi_t+\frac{(\varphi_x+\varphi_y+1)^2}{2}=0 \\
        \displaystyle    \varphi(x,y,0)=-\cos\left (\frac{\pi(x+y)}{2} \right)   \\
         \end{array}
   \right.
\end{equation}
with periodic boundary condition on the domain $[-2, 2]^2$.

In this example, we test the performance of our method on unstructured meshes. A sample mesh used with characteristic  length $h=1/4$ is given in Figure \ref{fig:mesh}.
At $t=0.5/\pi^2$, the solution is still 
smooth. Numerical errors and order of accuracy using $P^2$ polynomials are listed in 
Table \ref{table 2dburgtr}, demonstrating the optimal order
of accuracy.  At $t=1.5/\pi^2$, the solution is no longer smooth.  Our scheme could capture the viscosity solution as shown in Figure \ref{fig 2dburgtr}.

\begin{figure}[h!]

\centering
    \includegraphics[width=0.6\textwidth]{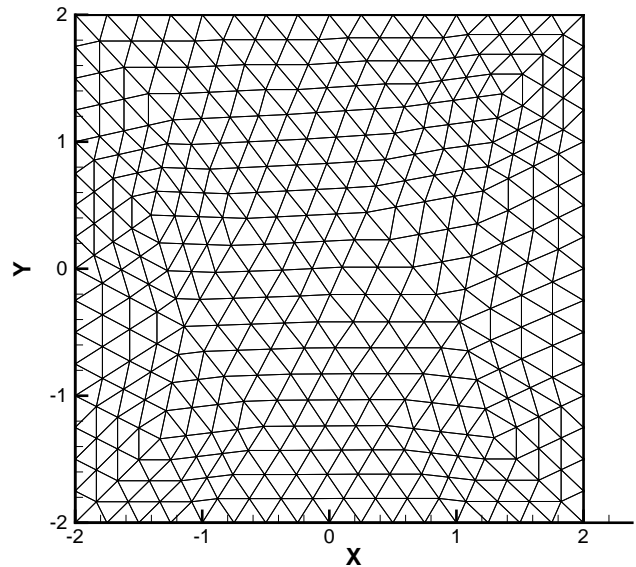}
 \caption{Examples \ref{2dburgtr} and \ref{2dcostr}. The unstructured mesh used with characteristic length $h=1/4$. }
  \label{fig:mesh}
\end{figure}

\begin{table}[ht]
\caption {Errors and numerical orders of accuracy for 
Example \ref{2dburgtr} when using $P^2$ polynomials and third order Runge-Kutta 
 time discretization on triangular meshes with characteristic length $h$. Penalty constant $C=0.25$.
Final time $t = 0.5/\pi^2$. $CFL=0.1$. }
\label{table 2dburgtr}
\vspace{2 mm}
\centering
\begin{tabular}{c c c c c c c}
\hline
$h$   & $L^1$ error & order & $L^2$ error & order & $L^\infty$ error & order  \\
\hline
 1 & 1.36E-02 &      & 2.31E-02 &       & 2.22E-01 &   \\
 1/2 & 1.77E-03 &  2.93 & 3.23E-03 &2.84 & 5.14E-02 &2.11\\
1/4 & 2.25E-04 & 2.98 &4.50E-04 &  2.84 & 8.95E-03  &  2.52  \\
1/8 & 2.74E-05 &3.04 & 5.82E-05 &  2.95 & 1.30E-03& 2.78  \\
1/16 & 3.40E-06 &3.01 & 7.53E-06 &  2.95 & 1.84E-04& 2.83  \\
\hline
\end{tabular}
\end{table}

\begin{figure}[h!]

\begin{minipage}{.5\textwidth}
  \centering
    \includegraphics[width=\linewidth]{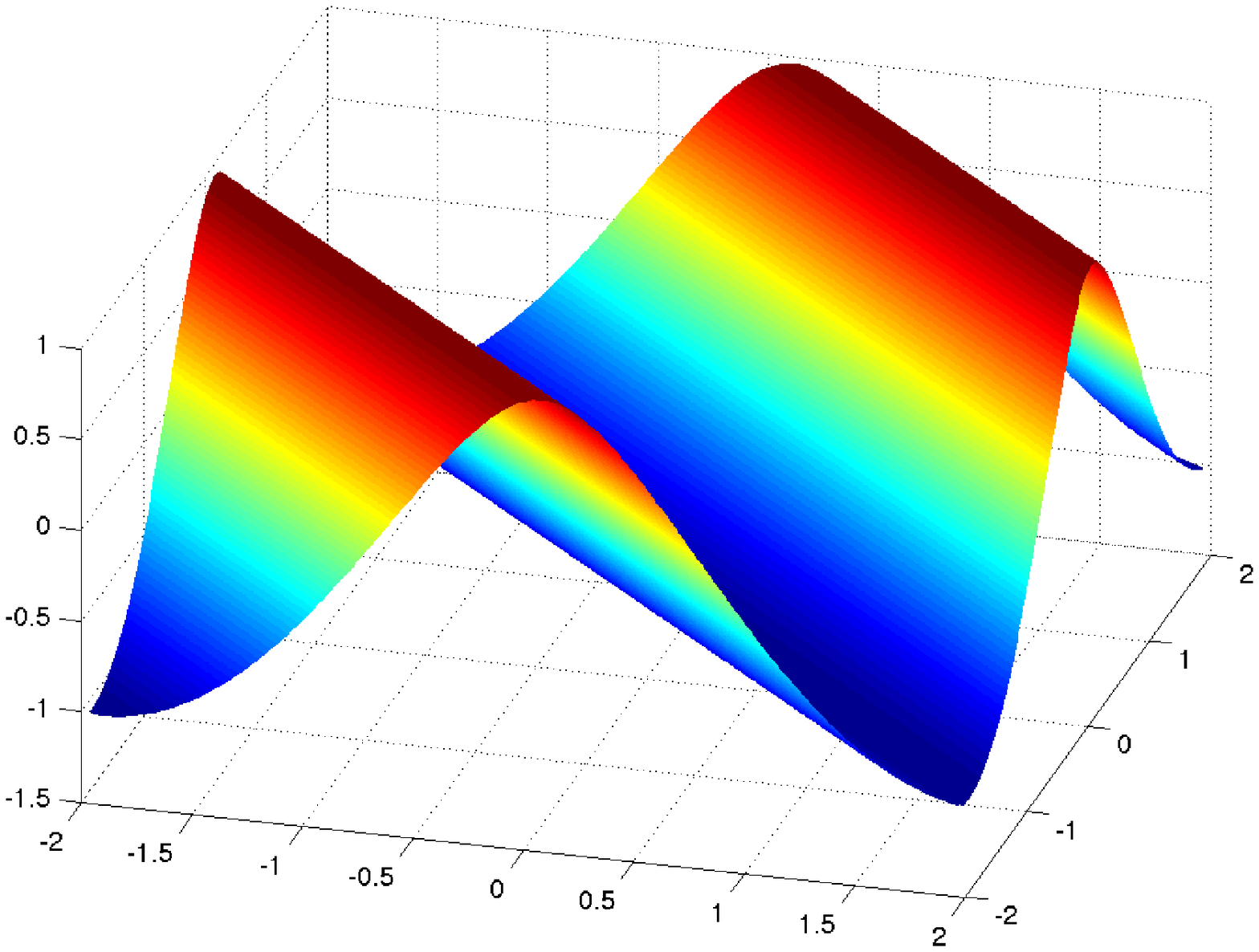}
\end{minipage}%
\begin{minipage}{.5\textwidth}
 \centering
  \includegraphics[width=\linewidth]{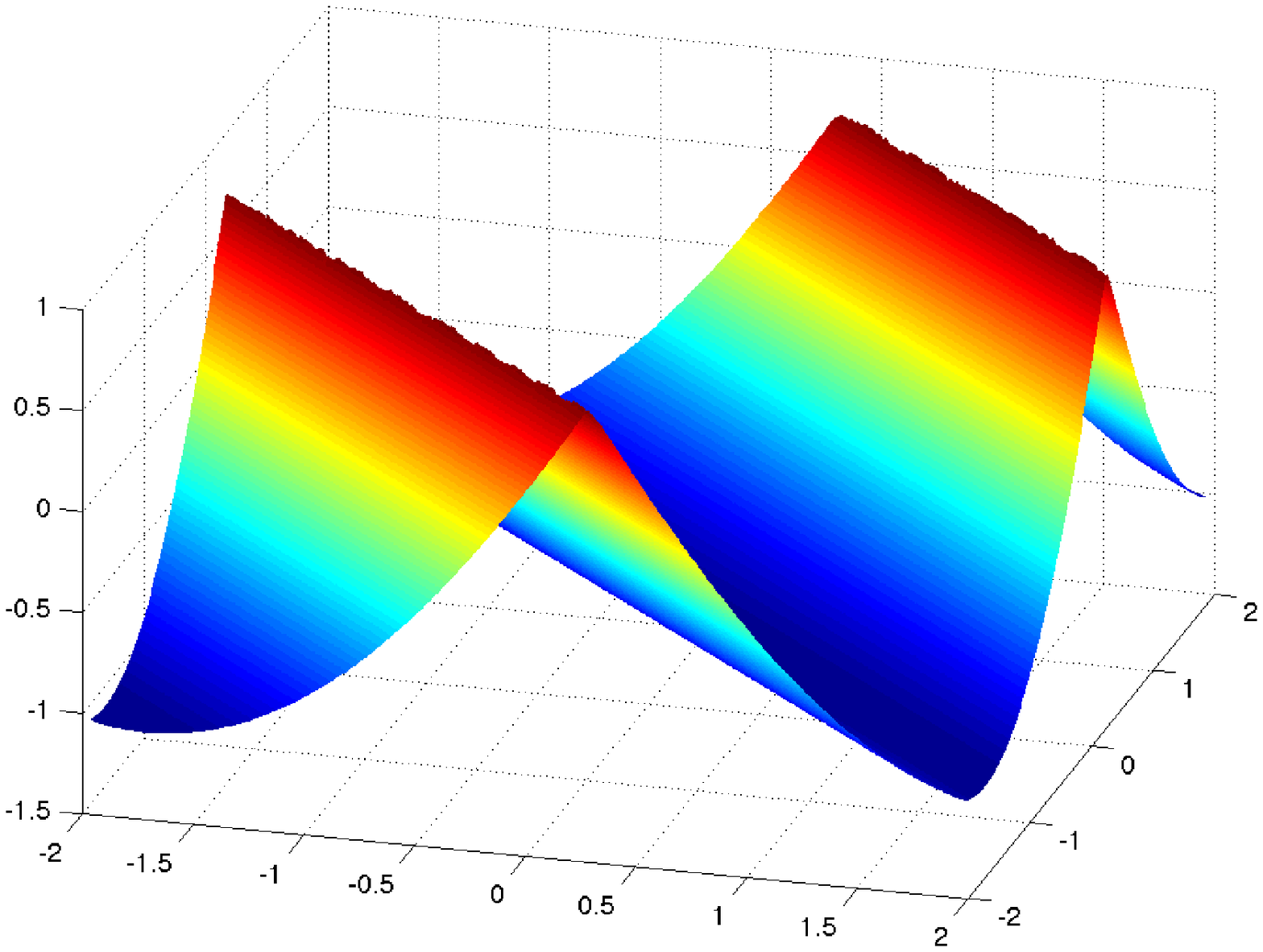}
\end{minipage}
 \caption{Example  \ref{2dburgtr}.   $CFL=0.1$.   $P^2$ polynomials. Triangular mesh with characteristic length $1/8$. 2816 elements. Penalty constant $C=0.25$. Left: $t=0.5/\pi^2$; right: $t=1.5/\pi^2$. }
  \label{fig 2dburgtr}
\end{figure}

\end{exa}

\begin{exa}
\rm 
\label{Li4.11} Two-dimensional nonlinear equation \cite{li2010central}
\begin{equation}
\left\{\begin{array} {l}
        \displaystyle    \varphi_t+\varphi_x\varphi_y=0 \\
        \displaystyle    \varphi(x,y,0)=\sin(x)+\cos(y)   \\
         \end{array}
   \right.
\end{equation}
with periodic boundary condition on the domain $[-\pi, \pi]^2$.

%Different from the previous two-dimensional problems, this example is genuinely two-dimensional. 
At $t=0.8$, the solution is still 
smooth, as shown in the left figure of Figure \ref{fig Li4.11}. Numerical errors and order of accuracy using $P^2$ polynomials are listed in 
Table \ref{table Li4.11}, demonstrating the optimal order
of accuracy.  At $t=1.5$, singular features would form in the solution, as shown in the right figure of Figure \ref{fig Li4.11}.  

%Table
\begin{table}[ht]
\caption {Errors and numerical orders of accuracy for 
Example \ref{Li4.11} when using $P^2$ polynomials and third order Runge-Kutta 
 time discretization on a uniform mesh of $N \times N$ cells. Penalty constant $C=0.25$.
Final time $t=0.8$. $CFL=0.1$. }
\vspace{2 mm}
\centering
\begin{tabular}{c c c c c c c}
\hline
N   & $L^1$ error & order & $L^2$ error & order & $L^\infty$ error & order  \\
\hline
 10 & 2.22E-03 &      & 3.95E-03 &       & 4.78E-02 &   \\
 20 & 2.75E-04 & 2.98 & 4.50E-04 & 2.98 & 7.79E-03 & 2.62\\
40 & 3.70E-05 & 2.89 & 7.33E-05  & 2.77 &1.50E-03  &  2.38  \\
80 & 4.80E-06 &2.95 & 9.83E-06 &  2.90 & 2.40E-04&  2.64  \\
\hline
\end{tabular}
\label{table Li4.11}.
\end{table}

% graph

\begin{figure}
\centering
\begin{minipage}{.5\textwidth}
  \centering
  \includegraphics[width=\linewidth]{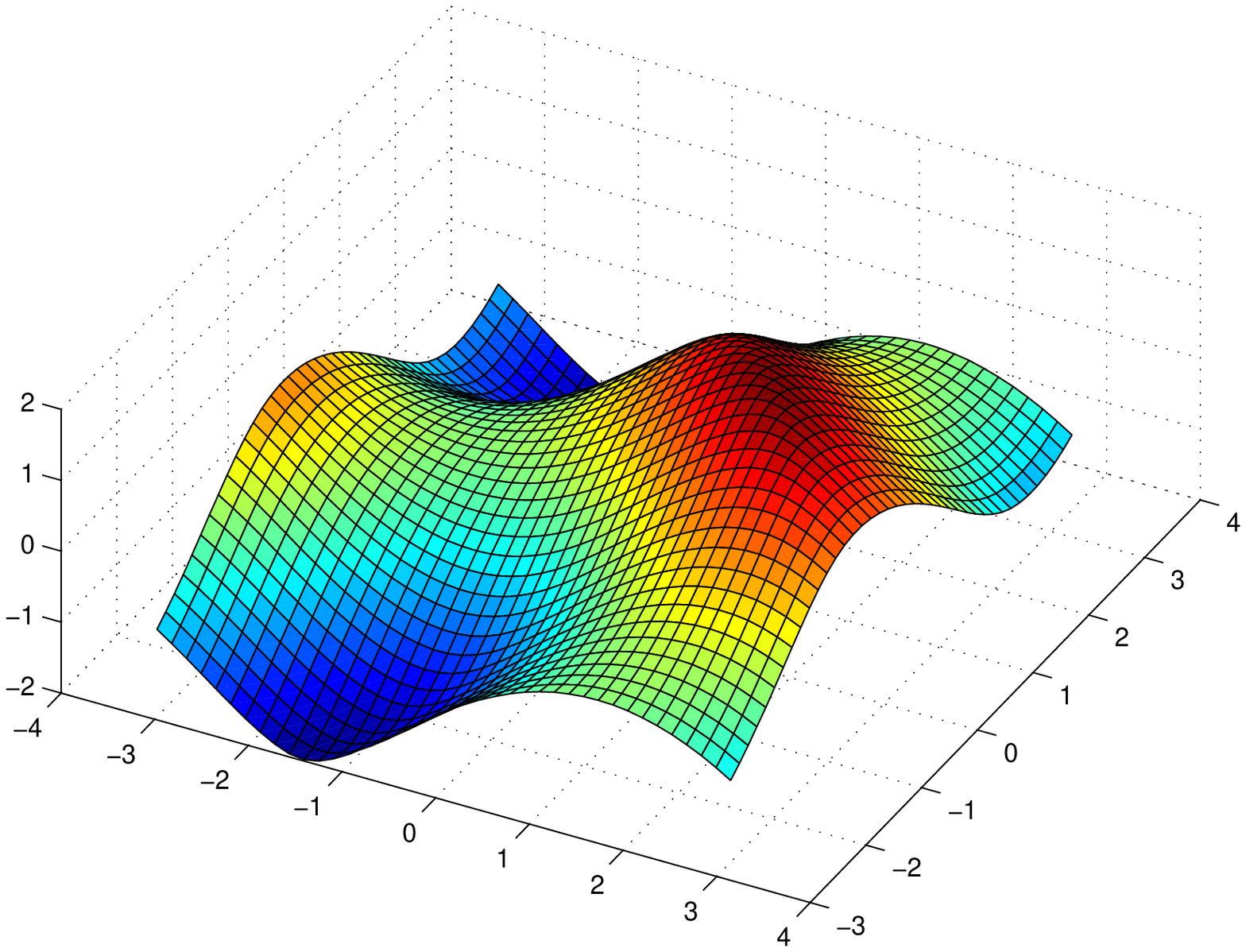}
  \label{Example 3.4}
\end{minipage}%
\begin{minipage}{.5\textwidth}
 \centering
  \includegraphics[width=\linewidth]{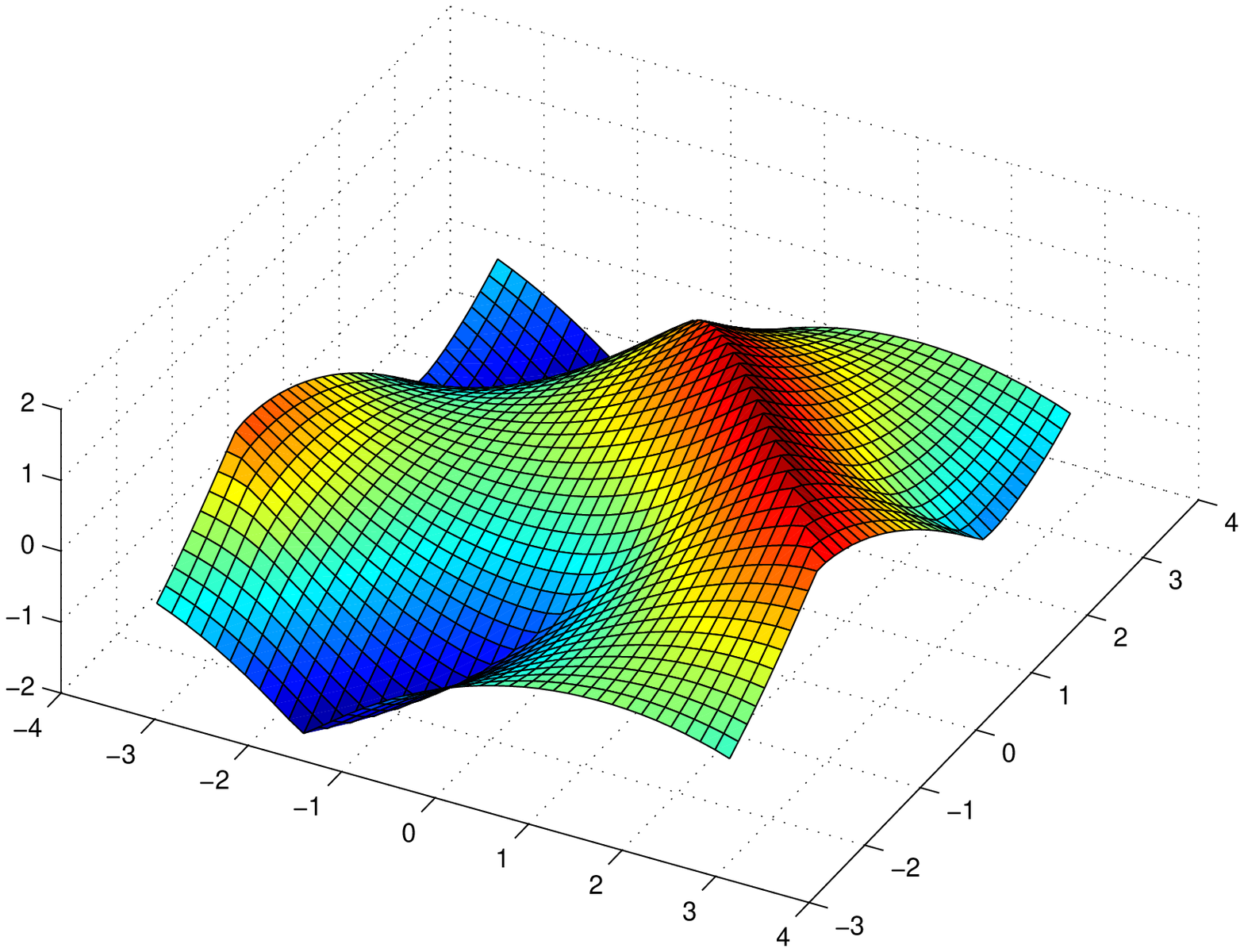}
  \label{Example 3.2}
\end{minipage}
\caption{Example \ref{Li4.11}. $CFL=0.1$.   $P^2$ polynomials on a $80\times 80$ uniform mesh. Penalty constant: $C=0.25$.    Left: $t=0.8$; right: $t=1.5$.  }
\label{fig Li4.11}
\end{figure}
\end{exa}

% Li Example 4.12 
\begin{exa}
\rm An example related to controlling optimal cost determination  \cite{Osher_1991_SIAM_NonOscill} \label{cost}
\begin{equation}
\left\{\begin{array} {l}
        \displaystyle    \varphi_t+\sin(y)\varphi_x+(\sin(x)+\sign(\varphi_y))\varphi_y-\frac{1}{2}\sin ^2(y)+\cos(x)-1=0 \\
        \displaystyle    \varphi(x,y,0)=0  \\
         \end{array}
   \right.
\end{equation}
\rm with periodic boundary condition on the domain $[-\pi, \pi]^2$.

\rm  The Hamiltonian is not smooth in this example. Our scheme can   capture the features of the viscosity solution well. The numerical solution (left) and the optimal control term $\sign(\varphi_y)$ (right) at $t=1$ are shown in Figure \ref{fig cost}.  

% graph

\begin{figure}
\centering
\begin{minipage}{.5\textwidth}
  \centering
  \includegraphics[width=\linewidth]{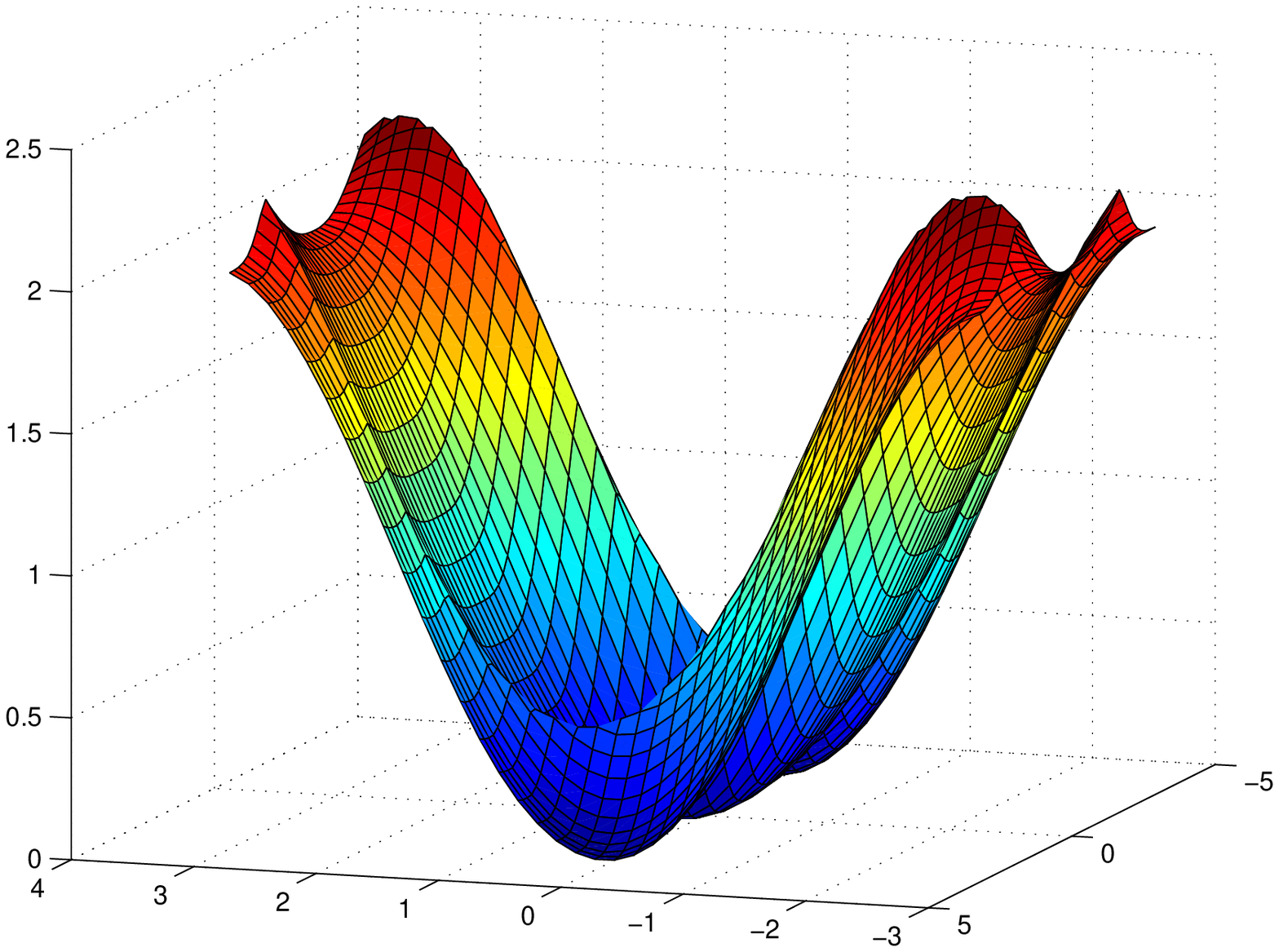}
  \label{Example 3.13}
\end{minipage}%
\begin{minipage}{.5\textwidth}
 \centering
  \includegraphics[width=\linewidth]{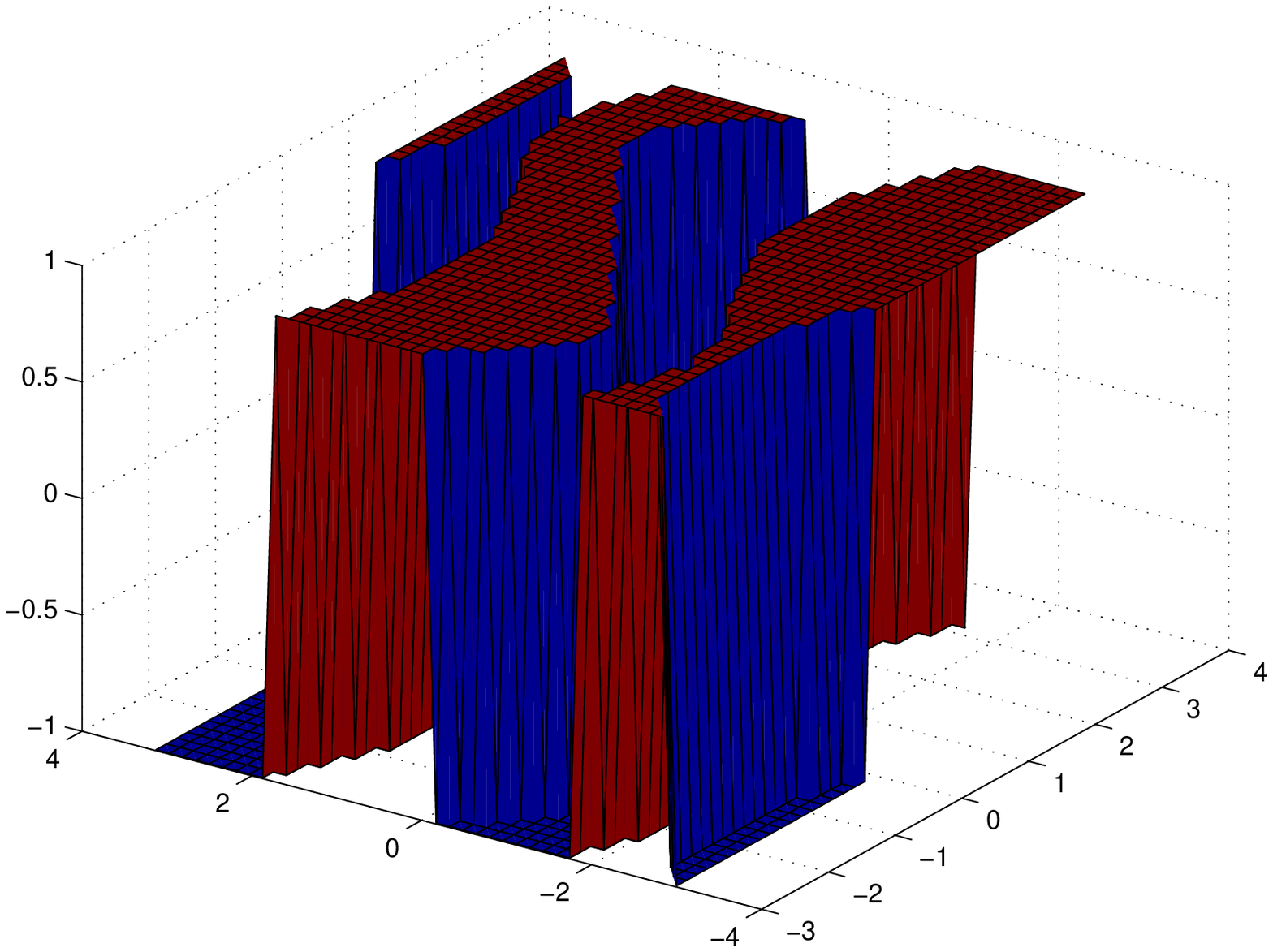}
  \label{Example 3.13}
\end{minipage}
\caption{Example \ref{cost}. $t=1$. $CFL=0.1$.  $P^2$ polynomials on a $40\times 40$ uniform mesh. Penalty constant: $C=0.25$.  Left: the numerical solution; right: $\sign(\varphi_y)$.  }
\label{fig cost}
\end{figure}
\end{exa}

\begin{exa}
\label{2dcostr}
\rm  Two-dimensional equation with a nonconvex Hamiltonian 
\begin{equation}
\label{eqn 2dcostr} \left\{\begin{array} {l}
        \displaystyle    \varphi_t-\cos(\varphi_x+\varphi_y+1)=0 \\
        \displaystyle    \varphi(x,y,0)=-\cos(\frac{\pi }{2} (x+y))   \\
         \end{array}
   \right.
\end{equation}
with periodic boundary condition on the domain $[-2, 2]^2$.

We use the same unstructured mesh as in Example \ref{2dburgtr}, see for example Figure \ref{fig:mesh}. At $t = 0.5/\pi^2$, the solution is still smooth, see Table \ref{table2dcostr} for numerical errors and order of accuracy using $P^2$ polynomials.  At $t = 1.5/\pi^2$, singular features would develop in the solution, as shown in Figure \ref{fig 2dcostr}.  

\begin{table}[ht]
\caption {Errors and numerical orders of accuracy for 
Example \ref{2dcostr} when using $P^2$ polynomials and third order Runge-Kutta 
time discretization on triangular meshes with characteristic length $h$. Penalty constant $C=0.25$.
Final time $t = 0.5/\pi^2$. $CFL=0.1$.  }
\label{table2dcostr}
\vspace{2 mm}
\centering
\begin{tabular}{c c c c c c c}
\hline
$h$   & $L^2$ error & order & $L^2$ error & order & $L^\infty$ error & order  \\
\hline
 1 & 1.05E-02 &      & 1.65E-02 &       & 1.48E-01 &   \\
 1/2 & 1.59E-03 &  2.71 & 2.49E-03 &2.73 & 3.11E-02 &2.25\\
1/4 & 2.42E-04 & 2.71 &4.02E-04 &  2.63 & 6.35E-03  &  2.29  \\
1/8 & 3.28E-05 &2.89 & 5.84E-05 &  2.78 & 1.03E-03& 2.62  \\
1/16 & 3.96E-06 &3.05 & 7.45E-06 &  2.97 & 1.72E-04& 2.58  \\
\hline
\end{tabular}
\end{table}

\begin{figure}[h!]
\begin{minipage}{.46\textwidth}
  \centering
\includegraphics[width=\linewidth]{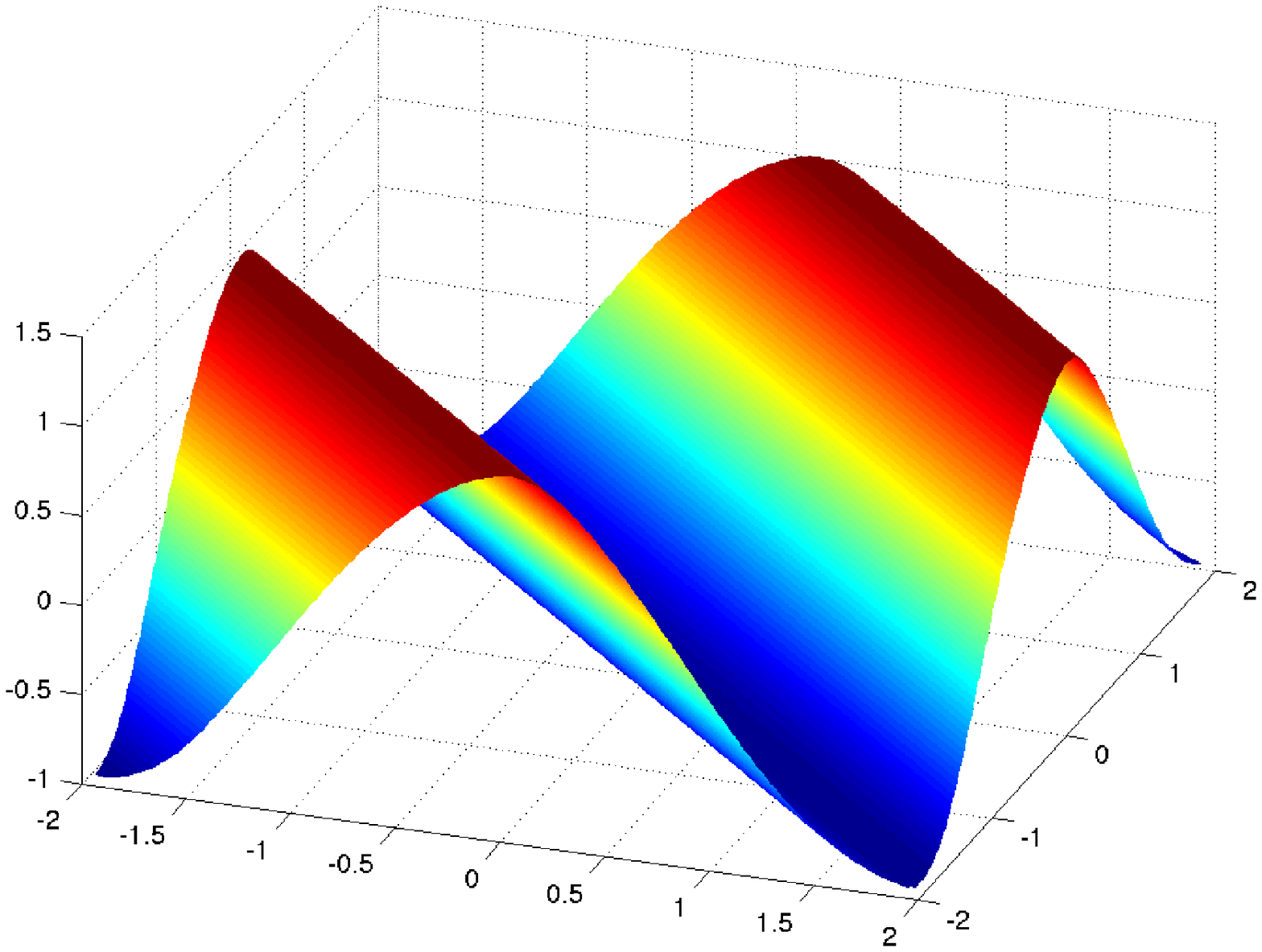}
\end{minipage}%
\begin{minipage}{.47\textwidth}
 \centering
\includegraphics[width=\linewidth]{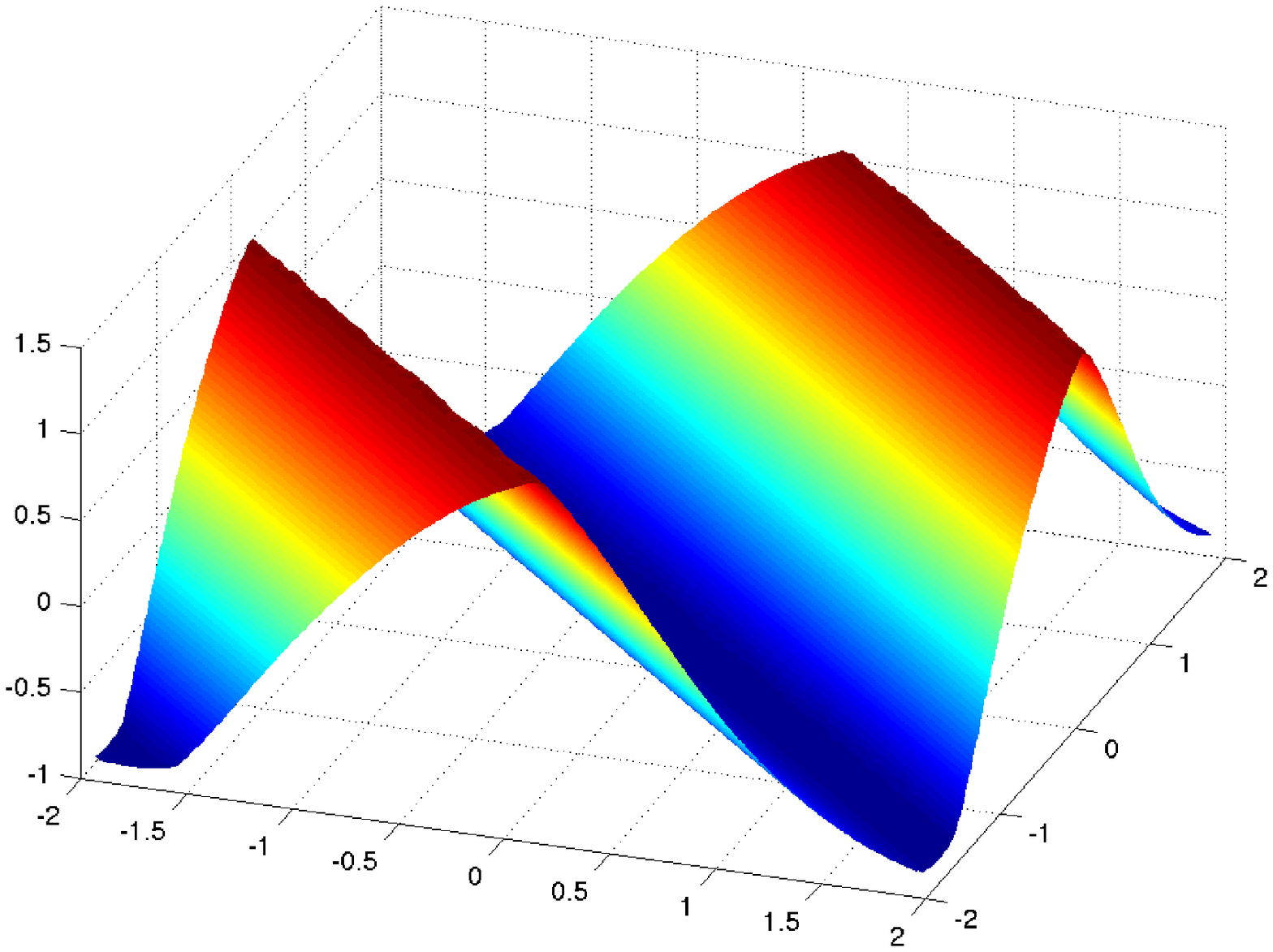}
\end{minipage}
\caption{Example  \ref{2dcostr}.  $CFL=0.1$.   $P^2$ polynomials. Triangular mesh with characteristic length $1/8$. 2816 elements. Penalty constant: $C=0.25$.  Left: $t=0.5/\pi^2$. Right: $t=1.5/\pi^2$.}
\label{fig 2dcostr}
\end{figure}
\end{exa}

%Yan Jue Example 4.9
\begin{exa}
\label{2dReimann}
\rm  Two-dimensional Riemann problem 
\begin{equation}
\label{2d.burg} \left\{\begin{array} {l}
        \displaystyle    \varphi_t+\sin(\varphi_x+\varphi_y)=0 \\
        \displaystyle    \varphi(x,y,0)=\pi (|y|-|x|)   \\
         \end{array}
   \right.
\end{equation}
on the domain $[-1, 1]^2$.

%\rm If $N$ is even, even with Moment limiter, it doesn't change, just the initial condition. If $N$ is odd, for example, $N=41$, then only for $P^2$.
\rm Similar to \cite{Hu_1999_SIAM_DG_FEM, yan2011local},  a nonlinear limiter is needed for convergence in this example. We use the moment limiter \cite{krivodonova2007limiters} and the numerical solution obtained by $P^2$ polynomial at $t=1$ is provided in Figure \ref{fig 2dReimann}.
%The moment limiter works in a similar way as Minmod limiter, it limits the derivative of order $i$ in a given cell using the derivatives of order $i-1$ in neighboring cells. 
% graph

\begin{figure}[h!]
\centering
    \includegraphics[width=0.6\textwidth]{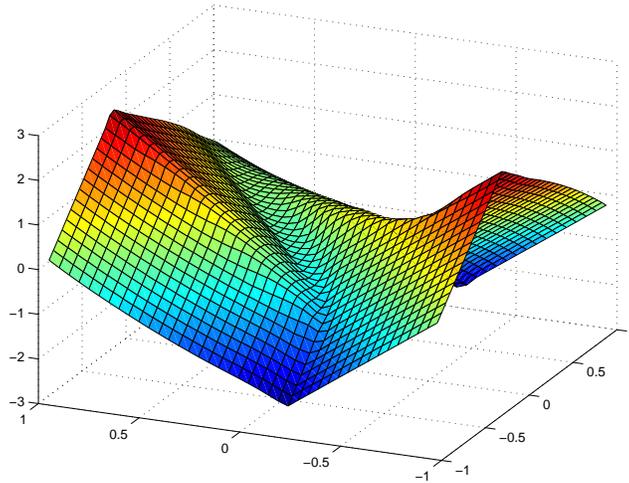}
 \caption{Example \ref{2dReimann}. $t = 1$. $CFL=0.1$.  $P^2$ polynomials on a $41\times 41$ uniform mesh. Penalty constant: $C=0.25$.  }
  \label{fig 2dReimann}
\end{figure}
\end{exa}

\begin{exa}
\label{2dprop}
\rm  The problem of a propagating surface  
\begin{equation}
\label{2d.prop} \left\{\begin{array} {l}
        \displaystyle    \varphi_t-\sqrt{\varphi_x^2+\varphi_y^2+1}=0 \\
        \displaystyle    \varphi(x,y,0)=1-\frac{1}{4}(\cos(2\pi x)-1)( \cos(2 \pi y)-1)  \\
         \end{array}
   \right.
\end{equation}
with periodic boundary condition on the domain $[0,1]^2$. This is a special case of the example used in \cite{osher1988fronts}, and  is also the two-dimensional Eikonal equation arising from geometric optics \cite{jin1998numerical}.
We use an unstructured mesh shown in Figure \ref{fig meshprop} with refinements near the center of the domain where the solution develops singularity. The numerical solutions at different times are displayed in Figure \ref{fig prop}. Notice that the solution at $t=0$ is shifted downward by 0.35 to show the detail of the solutions at later times.

\begin{figure}[h!]
\centering
    \includegraphics[width=0.6\textwidth]{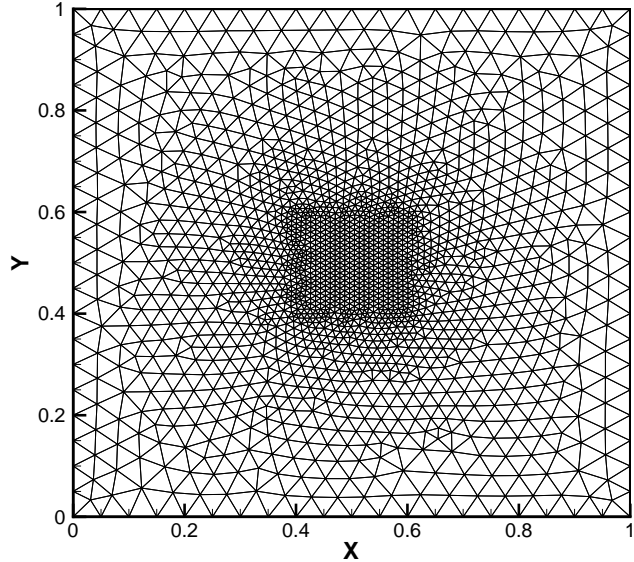}
     \caption{Example \ref{2dprop}. The unstructured mesh used in the computation. Number of elements: 3480. }
 \label{fig meshprop}
\end{figure}

\begin{figure}[h!]
\centering
    \includegraphics[width=0.7\textwidth]{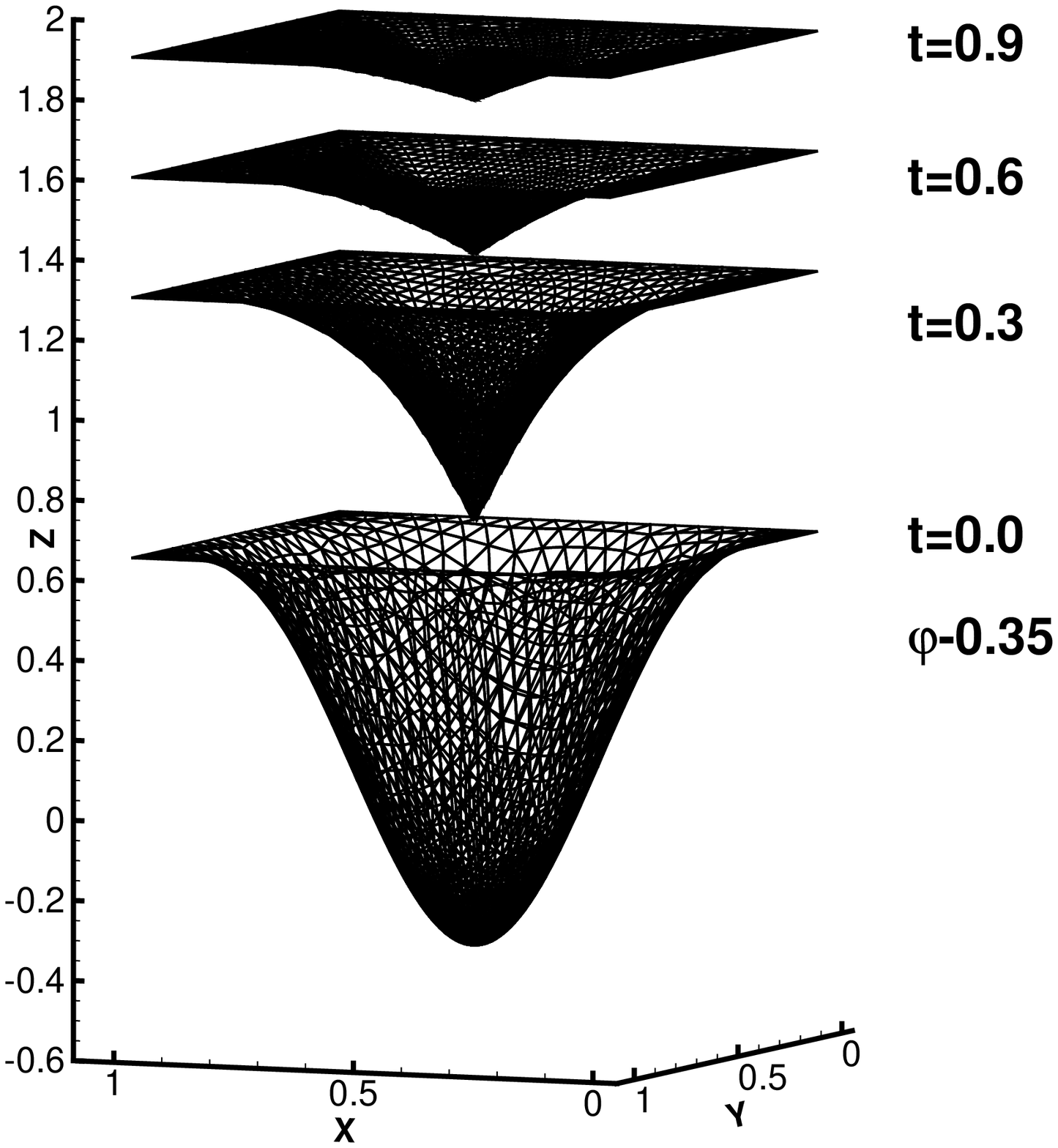}
\caption{Example \ref{2dprop}. $CFL=0.1$.   $P^2$ polynomials. Penalty constant: $C=0.25$. The numerical solution at the indicated times. }
 \label{fig prop}
\end{figure}

\end{exa}

\begin{exa}
\label{2ddisk}
\rm  The problem of a propagating surface on the unit disk $\{(x,y):x^2+y^2 \leq1\}$
\begin{equation}
\label{2d.disk} \left\{\begin{array} {l}
        \displaystyle    \varphi_t-\sqrt{\varphi_x^2+\varphi_y^2+1}=0 \\
        \displaystyle    \varphi(x,y,0)=-\sin(\frac{\pi (x^2+y^2)}{2}).   \\
         \end{array}
   \right.
\end{equation}
We use an unstructured mesh as depicted in Figure \ref{fig meshdisk} with refinements near the origin where solution develops singularity. The numerical solutions at different times are displayed in Figure \ref{fig disk}. Notice that the solution at $t=0$ is shifted downward by 0.2 to show the detail of the solutions at later times.

\begin{figure}[h!]
\centering
    \includegraphics[width=0.6\textwidth]{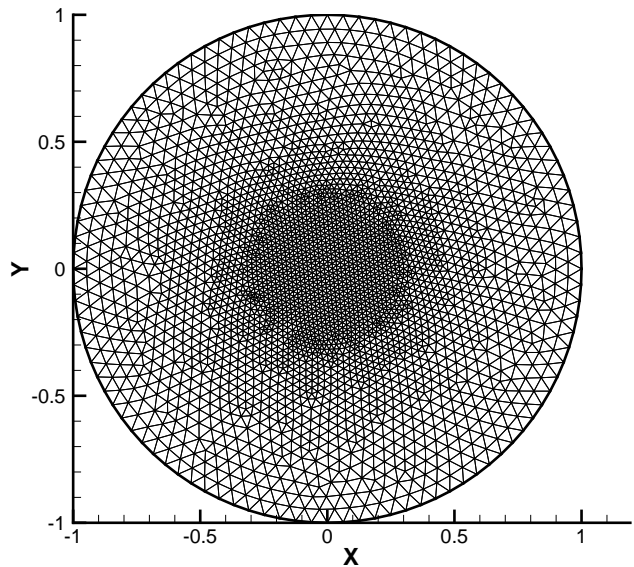}
 \caption{Example \ref{2ddisk}. The unstructured mesh used in the computation. Number of elements: 5890. }
  \label{fig meshdisk}
\end{figure}

\begin{figure}[h!]
\centering
    \includegraphics[width=0.7\textwidth]{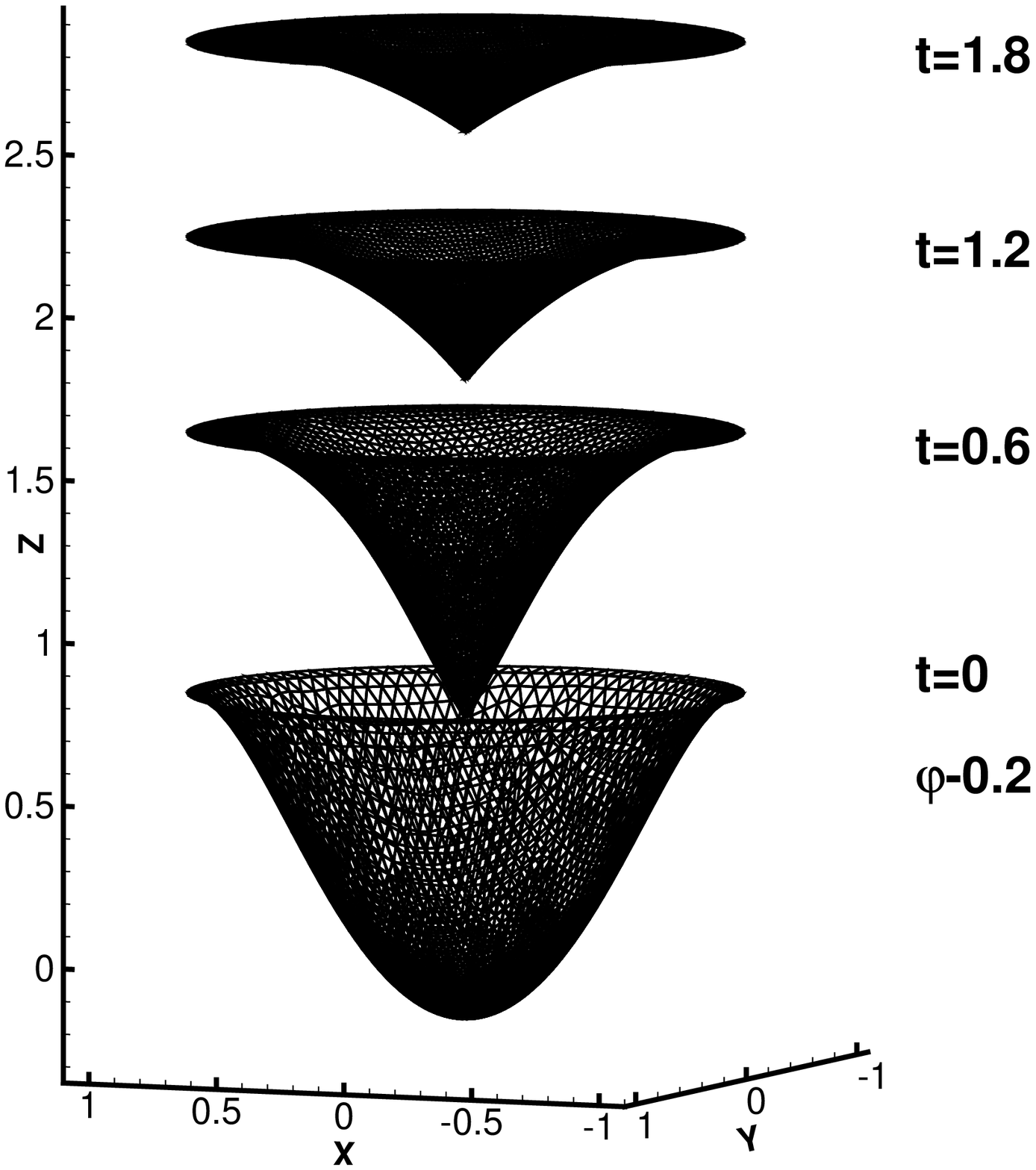}
 \caption{Example \ref{2ddisk}. $CFL=0.1$.   $P^2$ polynomials.  Penalty constant: $C=0.25$. The numerical solution at the indicated times. }
  \label{fig disk}
\end{figure}
\end{exa}

\section{Concluding Remarks}
\label{sec:conclusion}

In this paper, we propose a new DG method for directly solving the HJ equation. The scheme is direct and robust, and is  demonstrated to work on unstructured meshes even with nonconvex Hamiltonians. The theoretical aspects of  this method are subjects of future study.

\bibliographystyle{abbrv}
\bibliography{mypapers,ref_cheng}

\end{document}